\documentclass[twoside,11pt]{article}
\usepackage{jair}
\usepackage{theapa}
\usepackage{rawfonts}
\usepackage{style/my_packages}
\usepackage{style/my_commands}
\usepackage{style/my_operators}

\newcommand{\isdraft}{0}

\jairheading{1}{1993}{1-15}{6/91}{9/91}
\ShortHeadings{Distributed Bayesian: a continuous Distributed Constraint Optimization Problem solver}
{Fransman, Sijs, Dol, Theunissen \& De Schutter}
\firstpageno{0}

\begin{document}

\title{Distributed Bayesian: a continuous Distributed Constraint Optimization Problem solver}

\author{\name Jeroen Fransman \email j.e.fransman@tudelft.nl \\
       \addr Delft Center for Systems and Control
       (DCSC), Delft University of
       Technology
       \AND
       \name Joris Sijs \\ 
       \name Henry Dol \\ 
       \addr Netherlands Organisation for Applied
       Scientific Research (TNO)
       \AND
       \name Erik Theunissen \\
       \addr Netherlands Defence Academy
       (NLDA)
       \AND
       \name Bart De Schutter \\ 
       \addr Delft Center for Systems and Control
       (DCSC), Delft University of
       Technology
}

\maketitle

\begin{abstract}
    In this work, the novel \DBayLong algorithm is presented for solving multi-agent problems within the continuous \DCOPLong framework.
    This framework extends the classical DCOP framework towards utility functions with continuous domains.
    Traditional DCOP solvers discretize the continuous domains, which increases the problem size exponentially.
    \DBay overcomes this problem by utilizing \BO for the adaptive sampling of variables to avoid discretization entirely.
    We theoretically show that \DBay converges to the global optimum of the DCOP for Lipschitz continuous utility functions.
    The performance of the algorithm is evaluated empirically based on the sample efficiency.
    The proposed algorithm is compared to a centralized approach with equidistant discretization of the continuous domains for the sensor coordination problem.
    We find that our algorithm generates better solutions while requiring less samples.
\end{abstract}

\section{Introduction}
\label{section:introduction}

\structure{CONTEXT: Real world problems such as UAV: distributed, communication, limited processing power and memory}
Many real-world problems can be modeled as multi-agent problems in which agents need to assign values to their local variables to optimize a global objective based on utility.
Examples include scheduling \cite{Sato2015}, mobile sensor coordination \cite{Zivan2015}, hierarchical task network mapping \cite{Sultanik2007}, and cooperative search \cite{Acevedo2013}.
Even though numerous algorithms exist that solve these problems, applying them in practice is often problematic, as complications arise from limitations in communication, computation, and/or memory.

\structure{These problems are well suited to be modeled within the DCOP framework}
The \DCOPLong framework is well suited to model the above-mentioned problems (as detailed in \citeA{Meisels2008,Modi2005,Petcu2005b,Gershman2009,Yeoh2012}).
Within the DCOP framework a problem is defined based on variables and on utility functions that are aggregated into an objective function.
Additionally, agents assign values to all the variables that are allocated to them.
Agents are considered neighbors if their variables are arguments of the same utility function.
Neighbors then cooperatively optimize their utility functions through the exchange of messages.
Within the DCOP framework the variables are constrained by their domains.
In other words, a domain defines all possible value assignments of a variable.
This explicit definition of the domains of the variables is a major benefit for real-world problems that are (input) constrained.
Many solvers developed for DCOP assume that these domains are finite and discrete, while real-world problems are typically characterized by finite continuous domains.
A common approach for DCOP solvers is to use equidistant discretization, such as using a grid overlay to define all possible positions of an agent in an area.
When using equidistant sampling to discretize a continuous domain, the quality of the solution will depend on the distance between the values, where a smaller distance will allow for a better solution.
However, when continuous domains are discretized their cardinality will grow.
The increase in cardinality will exponentially increase the search space.
From the overview articles of \citeA{Leite2014} and \citeA{fioretto2019distributedreview}, it is clear that the cardinality of the domains is a major restriction to current DCOP solvers.
Therefore, discretization can cause solving continuous DCOPs to be intractable for current (discrete) DCOP solvers despite a small number of variables.

\structure{Solution is to take relation into account}
The underlying reason for the increase in problem size is that current DCOP solvers (implicitly) consider all values within a domain as unrelated to each other.
Because of this assumption it is not possible to efficiently sample the search space.
In problems with continuous domains this assumption can be removed since the utility of values that are close is often similar.
By explicitly taking this relation into account, a DCOP can be solved using efficient optimization methods.

\structure{Of possible approaches to take this relation into account, BO = best}
Various types of optimization methods exist in the literature that take advantage of the relation between the value and the utility.
Examples include simulated annealing \cite{cerny1985}, genetic algorithms \cite{Goldberg1989}, and \BO \cite{Mockus1989}.
In this work, \BO will be used as it focusses on efficient sampling during optimization, thereby requiring relatively few samples to closely approach the optimum.

\structure{Contributions of work}
Overall, the contributions of this article are threefold.
Firstly, we introduce an efficient algorithm that uses methods found in \BO to solve DCOPs with continuous domains called \DBayLong.
Secondly, we provide a theoretical proof of the convergence of the proposed algorithm called \DBay to the global optimum of the DCOP for utility functions with known Lipschitz constants.
Lastly, simulation results are given for a sensor coordination problem to compare the sample efficiency of the solver to a centralized approach based on equidistant discretization of the continuous domains.

\structure{Summary of article}
The remainder of this paper is organized as follows.
Firstly, in \cref{section:dcop_solvers} relevant literature regarding DCOP solvers are discussed.
Background information about the DCOP framework and the Bayesian optimization algorithm is provided in Sections \ref{section:DCOP} and \ref{section:BayesianOptimization}, respectively.
Afterwards, we present the novel sampling-based DCOP solver called \DBay in \cref{section:distributed_bayesian}.
The theoretical properties of \DBay are analyzed in \cref{section:analysis}.
Evaluation of \DBay for a sensor coordination problem is included in \cref{section:simulation_results}.
Finally, \cref{section:conclusion} summarizes the results and defines future work.
\clearpage  
\section{DCOP solvers}
\label{section:dcop_solvers}

\structure{Overview of DCOP, originate from CSP, therefore discrete}
The DCOP framework originates from an extension and generalization of Constraint Satisfaction Problems (CSPs) \cite{Tsang1993} towards distributed optimization.
A solution for a CSP is defined as the assignment of all variables from (finite) discrete domains such that all constraints are satisfied.
The CSP framework has been extended from centralized optimization to agent-based distributed optimization in the work of \citeA{Yokoo1998}.
Within the Distributed-CSP framework, the variables are allocated to agents and the agents coordinate the assignments of the variables among each other.

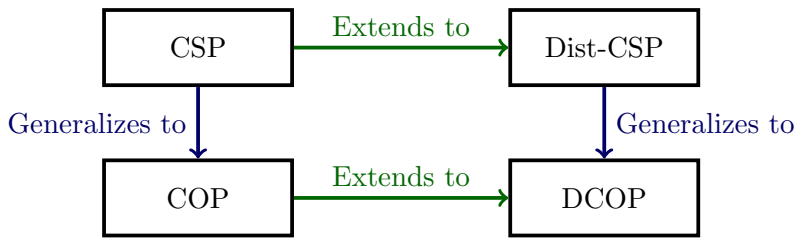
\begin{figure}[h]
    \centering
    \tikzstyle{classNode}=[
    rectangle,
    draw,
    line width=0.5mm,
    black,
    minimum height=1cm,
    minimum width=2.5cm,
]
\tikzstyle{classDistributedEdge}=[black!60!green,line width=0.5mm,->]
\tikzstyle{classUtilityEdge}=[black!60!blue,line width=0.5mm,->]

\def\classDistanceX{2.7}
\def\classDistanceY{1}

\begin{tikzpicture}[auto,-,thick]
	\node[classNode] (COP) at (-\classDistanceX,-\classDistanceY) {COP};
	\node[classNode] (CSP) at (-\classDistanceX,\classDistanceY) {CSP};
	\node[classNode] (DisCSP) at (\classDistanceX,\classDistanceY) {Dist-CSP};
	\node[classNode] (DCOP) at (\classDistanceX,-\classDistanceY) {DCOP};

    \path (CSP) edge[classDistributedEdge] node[classDistributedEdge] {Extends to} (DisCSP);
    \path (COP) edge[classDistributedEdge] node[classDistributedEdge] {Extends to} (DCOP);

    \path (CSP) edge[classUtilityEdge] node[classUtilityEdge,swap] {Generalizes to} (COP);
    \path (DisCSP) edge[classUtilityEdge] node[classUtilityEdge] {Generalizes to} (DCOP);
\end{tikzpicture}
    \caption{
        Graphical overview of the relations between the problem frameworks.
        Adapted from \citeA{fioretto2019distributedreview}.
        }
        \label{fig:COP_DCOP}
    \end{figure}

Additionally, CSP has been generalized into the Constraint Optimization Problem (COP) framework, where the constraints are replaced with utility functions.
These utility functions return a cost or reward based on the assignment of the variables.
Instead of constraint satisfaction, the goal of a COP is to find assignments that optimize an objective function.
Finally, the DCOP framework provides a unified framework that includes a large class of problems by combining the generalization and the extension.
A graphical overview of the relations between the problem frameworks can be seen in \cref{fig:COP_DCOP}.

However in the DCOP framework, the definition of the domains remains discrete, which limits its application to problems with continuous domains.
In the current paper, the (discrete) DCOP framework will be further extended to overcome this restriction in order to include real-world problems such as cooperative search.

\structure{How to extend/adapt the solvers for continious domains?}
As noted by \citeA{Modi2005}, optimally solving a discrete DCOP is NP-hard with respect to the number of variables and the cardinality of their domains.
For this reason, complete (optimal) DCOP solvers are often not used in practice.
In the literature a diverse range of incomplete (near-optimal) DCOP solvers exist that trade off solution quality against computational requirements.
Such solvers perform well for problems with domains a with low cardinality, such as graph coloring problems \cite{Modi2005} and meeting scheduling problems \cite{Jennings1995}.

Current DCOP solvers require the discretization of the continuous domains in the problem definition phase.
One typically discretizes the entire domain using a grid-based approach that converts the continuous domains into discrete domains.
This process can arbitrarily increase the cardinality of the domains, thereby rendering the use of current DCOP solvers intractable.
This raises the question whether there are DCOP solvers that can be extended to solve continuous DCOPs without the need for discretization.
\clearpage  
In the literature numerous solvers for DCOP with discrete domains have been proposed; for a detailed overview, the reader is referred to \citeA{Cerquides2014} and \citeA{Leite2014}.
In order to determine if DCOP solvers can be extended towards continuous domains, the taxonomy introduced by \citeA{Yeoh2008} is used.
The taxonomy defines three classes:
\begin{description}
    \item [Search-based solvers] perform a distributed search over the local search space of the agents.
    These solvers are based on centralized search techniques such as best-first and depth-first to \emph{reduce the search space} of the problem by exchanging messages between the agents.
    Examples are ADOPT \cite{Modi2005}, CoCoA \cite{Leeuwen2017}, AFB \cite{Gershman2009}, DSA \cite{Kirkpatrick1983}, and DBA \cite{Wittenburg2003}.
    \item [Inference-based solvers] communicate accumulated information among agents in order to \emph{reduce the problem size after every message} through dynamic programming methods.
    Well-known examples of this class of solvers are DPOP \cite{Petcu2005a}, the max-sum based algorithm \cite{Rogers2011}, and action GDL \cite{Vinyals2009}.
    \item [Sampling-based solvers] coordinate the sampling of the global search space guided by probabilistic measures.
    The probabilistic measures are calculated based on (all) preceding samples in order to \emph{balance exploration and exploitation} of the global search space.
    At the time of writing, two sample-based solvers are found in the literature: DUCT \cite{Ottens2018} and Distributed Gibbs \cite{Nguyen2019}.
\end{description}

In order to efficiently solve a continuous DCOP the relation between the value of the variables and the corresponding (global) utility values should to be taken into account explicitly.
Neither search-based nor inference-based solvers can be extended to take this property into account because both solver types only compare the utility of value assignments within a single iteration.
On the contrary, sampling-based solvers take previous iterations into account when selecting an additional sample.
This feature can be extended to take the relation between the samples into account during the sample selection process.
Which would allow for selecting a sample from a continuous domain directly.
Therefore, in this work, a \emph{novel sampling-based solver} will be introduced that removes the need for discretization entirely.

\section{Distributed Constraint Optimization Problems}
\label{section:DCOP}
A \DCOPLong is a problem in which an objective function needs to be optimized in a distributed manner through value assignments for all variables.
The objective function consists of the aggregate of utility functions, which define the utility value of the value assignments.
Agents are able to assign a value for variables that are allocated to them.
Furthermore, a variable can only be allocated to a single agent.
Typically, the number of variables is equal to the number of agents, i.e. every agent assigns a single variable.
The agents cooperate by sending messages to agents with whom they share a utility function.
A utility function is shared between agents if their variables are in the arguments of that function.
An important aspect of a DCOP is the definition of the domains of the variables.
A domain defines all possible values a variable can be assigned to.
In other words, the value assignments are restricted by the domains of the variables.

Following the notation of \citeA{fioretto2019distributedreview}, a continuous DCOP is defined by \mbox{\DCOPTupleLong} where
\begin{itemize}
    \item \agentSetLong is the set of agents,
        where \numberOfAgents is the number of agents.

    \item \variableSetLong is the set of variables, where $\numberOfVariables \geq \numberOfAgents$ is the number of variables.

    \item \domainSetLong is the set of domains of all variables, where $\domainSetIndexed \subseteq \realsSet$ is the (continuous) domain associated with variable \variableIndexed.
        The search space of the DCOP is defined by all possible combinations of all values within the domains as \searchSpaceLong, where $\prod$ is the set Cartesian product.
        The search space of a set of variables ($\variableSetFunction \subseteq \variableSet$) is defined as $\searchSpace_\variableSetFunction = \underset{\variableIndex ~:~ \variableIndexed \in \variableSetFunction}{\prod} \domainSetIndexed$.

        An assignment denotes the projection of variables onto their domain as \assignmentCompleteLong.
        In other words, for all $\variableIndexed \in \variableSet$ if \assignmentIndexed is defined, then $\assignmentIndexed \in \domainSetIndexed$.
        An assignment of a subset of variables is denoted by \assignmentPartialLong.

    \item \functionSetLong is the set of utility functions,
        where \numberOfFunctions is the number of utility functions.
        The scope of \utilityFunctionIndexed is denoted as \variableSetFunctionIndexedLong, where $\variableIndexed \in \variableSetFunctionIndexed$ when \variableIndexed is an argument of \utilityFunctionIndexed.
        The optimum of a utility function is denoted by \functionOutputOptimalLong with input \functionInputOptimalLong, where $\searchSpaceFunctionIndexed$ denotes the domain of the utility function.

    \item \variableMappingLong is a mapping from variables to agents.
        The agent to which variable \variableIndexed is allocated is denoted as \variableMappingAgentIndexed.
        A common assumption is that the number of agents is equal to the number of variables, such that $\agentIndexed = \variableMappingAgentIndexed$ for $\variableIndex = 1, \dots, \numberOfVariables$.
        Likewise, the set of agents associated with \utilityFunctionIndexed is denoted by \functionScopeLongIndexed.

    \item $\operator$ is an operator that combines all utility functions into the objective function.
        Common options are the \emph{summation} operator $\bracketBig{\sum(\cdot)}$ and the \emph{maximum} operator $\big(\max(\cdot) \big)$.
        The objective function is defined by \globalUtilityLong.
        The optimal value assignment is denoted by \assignmentOptimalLong.

\end{itemize}

\begin{figure}[h]
    \centering
    \begin{subfigure}[b]{0.3\textwidth}
        \tikzstyle{variableNode}=[rectangle,draw,line width=0.5mm,blue]
\tikzstyle{agentNode}=[circle,draw=black,line width=0.5mm,minimum size=1cm]
\tikzstyle{functionEdge}=[black!60!green,line width=0.5mm]

\begin{tikzpicture}[auto,-,thick]
	\node[agentNode,label=90:$a_1$] (a1) at (0,0) {};
	\node[agentNode,label=90:$a_2$] (a2) at (-1.5,-1) {};
	\node[agentNode,label=90:$a_3$] (a3) at (1.25,1.5) {};
	\node[agentNode,label=90:$a_4$] (a4) at (-1,2) {};
	\node[agentNode,label=90:$a_5$] (a5) at (2,-0.75) {};			

	\node[variableNode] (x1) at (a1) {$x_1$};
	\node[variableNode] (x2) at (a2) {$x_2$};		\node[variableNode] (x3) at (a3) {$x_3$};
	\node[variableNode] (x4) at (a4) {$x_4$};
	\node[variableNode] (x5) at (a5) {$x_5$};			

	\path (a1) edge[functionEdge] node[functionEdge] {$f_1$} (a2);
	\path (a1) edge[functionEdge] node[functionEdge,swap] {$f_2$} (a3);	
	\path (a3) edge[functionEdge] node[functionEdge,swap] {$f_3$} (a4);
	\path (a1) edge[functionEdge] node[functionEdge] {$f_4$} (a4);
	\path (a3) edge[functionEdge] node[functionEdge] {$f_5$} (a5);			
\end{tikzpicture}
        \caption{Constraint graph}
        \label{fig:dcop_constraint_graph}
    \end{subfigure}
    \qquad \qquad
    \begin{subfigure}[b]{0.3\textwidth}
        \tikzstyle{variableNode}=[rectangle,draw,line width=0.5mm,blue]
\tikzstyle{agentNode}=[circle,draw=black,line width=0.5mm,minimum size=1cm]
\tikzstyle{functionEdge}=[black!60!green,line width=0.5mm]

\begin{tikzpicture}[auto,-,thick]
	\node[agentNode,label=90:$a_1$] (a1) at (0,0) {};
	\node[agentNode,label=90:$a_2$] (a2) at (-1.5,-1.5) {};
	\node[agentNode,label=90:$a_3$] (a3) at (1.5,-1.5) {};
	\node[agentNode,label=95:$a_4$] (a4) at (0.25,-3) {};
	\node[agentNode,label=90:$a_5$] (a5) at (2.75,-3) {};			

	\node[variableNode] (x1) at (a1) {$x_1$};
	\node[variableNode] (x2) at (a2) {$x_2$};		\node[variableNode] (x3) at (a3) {$x_3$};
	\node[variableNode] (x4) at (a4) {$x_4$};
	\node[variableNode] (x5) at (a5) {$x_5$};			

	\path (a1) edge[functionEdge] node[functionEdge,swap] {$f_1$} (a2);
	\path (a1) edge[functionEdge] node[functionEdge] {$f_2$} (a3);	
	\path (a3) edge[functionEdge] node[functionEdge] {$f_3$} (a4);
	\path (a1) edge[functionEdge,dotted] node[functionEdge,swap,pos=.2] {$f_4$} (a4);
	\path (a3) edge[functionEdge] node[functionEdge,swap] {$f_5$} (a5);			
\end{tikzpicture}
        \caption{Pseudo-tree}
        \label{fig:dcop_pseudo_tree}
    \end{subfigure}
    \caption{
        Graphical example for two representations of the same DCOP.
        A node represents an agent $\agentIndexed$ and the edges indicate the utility functions $\utilityFunctionIndexed$ based on the variables $\variableIndexed$.
        In \cref{fig:dcop_constraint_graph} the nodes are unstructured, while in \cref{fig:dcop_pseudo_tree} the hierarchy is indicated through separate layers.
        The agent at the top ($\agent_\agentIndexRoot$) is the root of the tree, and the bottom agents ($\agent_2, \agent_4, \agent_5$) are the leaves.
        Note that the edge between $\agent_\agentIndexRoot$ and $\agent_4$ specifies a pseudo-parent/pseudo-child relation.
        }
        \label{fig:DCOP_representations}
\end{figure}
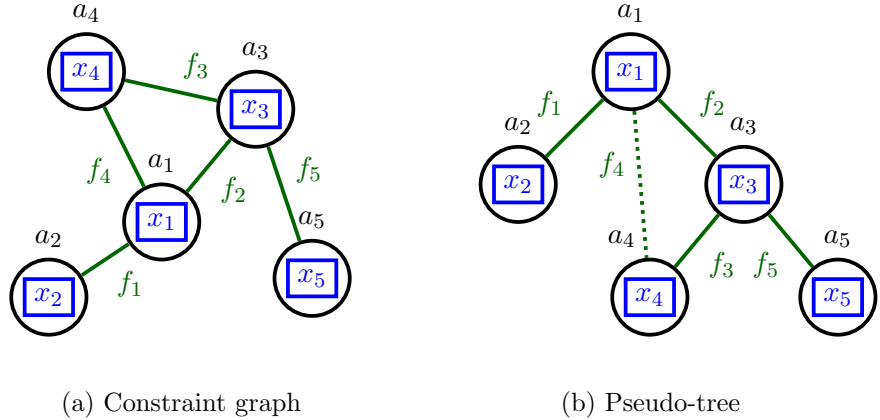

The relation between the variables and the utility functions is typically represented as a constraint graph.
In this representation the agents are defined as nodes and the edges implicitly represent the utility functions.
A constraint graph is often converted into a pseudo-tree \cite{Freuder1985} through a distributed depth-first-search procedure to partition the problem into independent subproblems based on the utility functions.
A pseudo-tree is a rooted spanning tree that introduces hierarchy to the agents where the subproblems are contained in separate \emph{branches}.
Additionally, all agents are assigned a single parent, which is an agent higher in the hierarchy.
The only exception is the agent on top of the hierarchy, which is called the root of the tree, this agent has no parent.
An agent can have multiple children and the agents without children correspond to the leaves of the tree.
In addition to parent/child relations, the pseudo-tree defines pseudo-parent/pseudo-child relations to indicate relations between agents over multiple hierarchy levels.
Typically, the pseudo-tree is used as a communication structure, where agents only communicate between parent and child.
In these cases, the pseudo relation allows for (indirect) interaction between pseudo-parent and pseudo-children.
A graphical example of the two DCOP representations is given in \cref{fig:DCOP_representations}.

\section{Sample selection for continuous DCOPs} 
\label{section:BayesianOptimization}

\structure{Which sample selection process to use?}
As mentioned in \cref{section:dcop_solvers}, in order to efficiently solve a continuous DCOP a sample-based solver is required.
To avoid discretization, a sample-based solver will be introduced such that the relation between the assignment and the utility will be taken into account within the sample selection process.

\structure{Overview of possibilities}
Within the literature several suitable methods exist for the sampling of continuous domains:
\begin{enumerate*}[label={\roman*)}]
    \item simulated annealing \cite{cerny1985}, which is a method that resembles the annealing of metals during cooling,
    \item genetic algorithms \cite{Goldberg1989}, which are a class of algorithms that resemble optimization through survival of the most \emph{fit} solutions, and
    \item Bayesian optimization \cite{Mockus1982}, which is a method to find the optimum of a function by using a probabilistic model of that function.
\end{enumerate*}
All methods mentioned above could potentially be used within a sample-based DCOP solver.
In this work \BO will be used as for certain probabilistic models convergence to the global optimum can be guaranteed as shown in \citeA[Theorem~6]{Vazquez2010}.
These results are not available for either simulated annealing or genetic algorithms.

Indeed, as noted by \citeA{Rios2013}, based on \citeA{Romeo1991}, although convergence results exist for simulated annealing no performance guarantees can be given for finite number of iterations.
Genetic algorithms, as summarized in \shortciteA{El-Mihoub2006}, are capable of converging to a region in which the global optimum exists; however they may require a large number of samples to converge to an optimum within this region \cite{de2005genetic}.
As highlighted by \citeA{Lee2018}, this optimum is not guaranteed to be the global optimum.

\structure{What is BO in detail?}
\BO is a method to find the global optimum of a function in a sample-efficient manner, i.e.\ it minimizes the number of required samples.
\BO consists of two elements: a \emph{probabilistic model} to approximate a (utility) function \utilityFunctionDot, and an \emph{acquisition function} \acquisitionFunctionDot to optimally select a new sample \nextSample, where \observationIndex denotes the sample index.
These two elements are discussed in more detail in the subsequent sections.
Based on samples of the function and the corresponding function values \functionOutputLongIndexed, the probabilistic model is used to estimate a mean function $\meanFunctionBase(\cdot)$ and the corresponding variance function $\varianceFunctionBase(\cdot)$.
Every input/output pair, \observationLongIndexed, is included in the ordered observations set \observationSetLong, where \numberOfObservations is the number of observations.
The observations are (re)ordered after every new observation, such that for scalar arguments \mbox{$\functionInputBase_1 \leq \functionInputBase_2 \leq \dots \leq \functionInputBase_\numberOfObservations$}.
The observations are used to update the probabilistic model, such that after every new observation the approximation is refined.
An overview of the Bayesian optimization algorithm is given in Algorithm \ref{algorithm:BO}.
\begin{algorithm}
    \SetKwInOut{Input}{Input}
    \SetKwInOut{Output}{Output}
    \Input{\utilityFunctionDot, \acquisitionFunctionDot, \functionInput}
    \Output{$\meanFunctionBase(\functionInput, \observationSet)$, $\varianceFunctionBase(\functionInput, \observationSet)$}
    \For{s = 1, 2, \dots}{
        \tcc{Select the next sample based on acquisition function}
        $\nextSample \coloneqq \argmax_{x} q(x | \mathcal{O}_{s-1})$\;
        \tcc{Sample the utility function}
        $y_s \coloneqq \utilityFunction(\nextSample)$\;
        \tcc{Augment (and reorder) the observation set}
        $\mathcal{O}_{s} \coloneqq \mathcal{O}_{s-1} \cup \{ (\nextSample, y_s) \}$\;
        \tcc{Calculate the mean function and the variance function}
        $\meanFunctionBase(\functionInput, \observationSet) =
        \expectation\left[ \utilityFunction(\functionInput) | \observationSet \right]
        $\;
        $\varianceFunctionBase(\functionInput, \observationSet) =
        \expectation\left[
            \left(
                [\utilityFunction( \functionInput) | \observationSet] - \meanFunctionBase(\functionInput, \observationSet)
            \right)^2
        \right]
        $\;
    }
    \caption{Bayesian optimization \cite{Mockus1982}}
    \label{algorithm:BO}
\end{algorithm}

\subsection{Probabilistic model}
\label{subsection:probabilistic_model}
The \GP{} is a widely used probabilistic model to represent acquired \mbox{knowledge} about a function.
More elaborate models exist, but these will often undo the computational benefit of the \GP model.
Using the \GP model, a function \utilityFunctionDot is modeled based on a prior mean function \gpMeanFunctionLong and a kernel \kernelLong.
The kernel represents the cross-correlation between the two variables \functionInputBase, $\functionInputBase^\prime$.
The prior mean function and the kernel contain all (prior) knowledge of \utilityFunctionDot.
Typically, the prior mean function is equal to the zero function ($\gpMeanFunction(\functionInputBase) = 0$ for all \functionInputBase) if no prior information about the function is available.
Therefore, the modelling of the function depends mostly on the choice of \kernel.

The \GP model is combined with the observations in order to construct the joint Gaussian distribution over the function.
From the joint Gaussian distribution, the posterior (distribution) can be found by using the Sherman-Morrison-Woodbury formula \cite{Sherman1950}:
\begin{align}
    \probability(\utilityFunction(\functionInputBase) | \observationSet) &= \NormalDistribution\bracket{\meanFunction, \varianceFunction }, \label{eq:posterior}
\intertext{where}
    \meanFunction &= \meanFunctionExtended \label{eq:mean_function} \\
    \varianceFunction &= \varianceFunctionExtended \label{eq:variance_function}
\end{align}
and \NormalDistribution denotes the normal distribution,
$\kernelGramian(\observationSet)$ is the Gramian matrix of the kernel, defined by \kernelGramianLong,
\kernelCrossCorrelationLong denotes the \emph{cross-correlation} vector between the observations and \functionInputBase,
and \observationOutputsLong denotes the observation value vector.
The (posterior) mean and variance function of the probabilistic model are denoted as \meanFunctionDot and \varianceFunctionDot, respectively.
Note that the posterior distribution contains the estimate of the function based on both the prior knowledge and the observations.

A wide range of kernels for Gaussian processes exist in the literature and the interested reader is referred to the work of \citeA{Duvenaud2011} for a overview on constructing kernels.
An important kernel property is universal approximation, which indicates the ability to estimate every continuous function up to a required resolution given a sufficient number of observations.
A kernel that possesses this property is called a universal kernel.
In the work of \citeA{Micchelli2006a}, the conditions for a kernel to be universal in terms of properties of its features are given.
Additionally, several examples of kernels are given that possess this property.
The most commonly used universal kernel is the \SEKernel.
A general description of the kernel and its properties is given by \citeA{Vert2004}.
A drawback of the \SEKernel is that it can result in over-smooth approximations.
For this reason, the \Matern kernel \cite{Minasny2005} is often used, since it can trade off differentiability and smoothness.
In practice, the choice for a kernel depends on the properties of the problem that needs to be solved.
Additionally, depending on the kernel one or more parameters need to be set.
All kernels have parameters that can be used to adjust their properties, such as smoothness and scaling.
If information about \utilityFunctionDot is available, this should be incorporated in the selection of the kernel and its parameters.
Typically, it is assumed that no information about \utilityFunctionDot is available, and then, as noted by \citeA{Rasmussen2004}, the selection of the parameters is non-trivial.
For this reason, the selection of parameters is often treated as a separate optimization problem \cite{MacKay1992}.
It is commonly solved by using the maximum likelihood criterion for which \ARD \cite{Mackay1994} is a widely used algorithm.

\subsection{Acquisition function}
\label{subsection:acquisition_function}
The selection of the next sample is the result of the optimization of an acquisition function \acquisitionFunctionDot, defined by
\begin{align*}
    \nextSample &= \argmax_{\functionInputBase} \acquisitionFunctionBase(\functionInputBase | \observationSet).
\end{align*}
The acquisition function depends on the posterior distribution in \cref{eq:posterior} and thereby on all previous observations.
Two commonly used acquisition functions are the \PI function \cite{Kushner1964} and the \EI function \cite{Mockus1978}.
The probability of improvement function considers the \emph{probability} of finding an observation the value of which is larger than the maximum observed value.
The maximum observed value is defined as
\begin{align*}
    \observationMaxSampleLong.
\end{align*}
The corresponding input is defined as \observationMaxSampleInputLong.
As noted by \citeA{brochu2010tutorial}, the probability of improvement function focusses solely on exploitation of already observed samples.
In order to balance exploration of the search space and exploitation of the observations, the expected improvement function will be used in this work.
The expected improvement function chooses the sample based on the expected \emph{value} of the next observation.
The expected improvement function can be written in closed form in terms of the mean and the variance function of the probabilistic model as
\begin{align}
    \acquisitionFunction &=
    \begin{cases}
        \acquisitionFunctionNonZeroCumulativeLong
        +
        \acquisitionFunctionNonZeroDensityLong
        &\text{if }
        \deviationFunction > 0\\
        0 &\text{if } \deviationFunction = 0
    \end{cases}
    \label{eq:acquisition}
    \\
    \acquisitionHelper &= \acquisitionHelperExtended
\end{align}
where $\Phi(\cdot)$ is the Gaussian cumulative distribution function, $\phi(\cdot)$ is the Gaussian probability density function, and \acquisitionFunctionParameter is a design parameter.
The design parameter can be used to tradeoff exploration and exploitation.
As noted by \citeA{Lizotte2012}, even for values as low as $\acquisitionFunctionParameter = 0$ will not result in a solely exploiting sampling behavior.
The interested reader is referred to \citeA{brochu2010tutorial} for a comparison of the two acquisition functions and more details.

\section{The Distributed Bayesian algorithm}
\label{section:distributed_bayesian}

\structure{Small introduction}
In the previous sections, background information has been given about the DCOP framework and \BO.
In this section, the novel sample-based DCOP solver \DBayLong is presented.
This solver is capable of directly solving continuous DCOPs without discretization of the domains.
\DBay uses \BO as the probabilistic measure to optimize the sample selection.
The overall procedure is similar to state-of-the-art sample-based solvers, e.g. DUCT \cite{Ottens2018} and \Gibbs \cite{Nguyen2019}.

\structure{Sample-based solvers in general}
Sample-based solvers coordinate the sampling of the global search space guided by probabilistic measures in order to balance exploration and exploitation.
The general outline of sample-based solvers is as follows.
Based on a pseudo-tree representation of the DCOP, the variables and utility functions are allocated to the agents.
Afterwards, two consecutive phases are iteratively repeated until termination conditions are satisfied.
The first phase, the sampling phase, is top-down and starts from the root agent.
The root starts the sampling phase by selecting a sample for all its variables.
A sample can be viewed as a temporary assignment of a variable.
The sample is sent as a \sampleMessageText message to all the children of the agent.
Upon receiving this message, an agent samples its own variables and adds these samples to the \sampleMessageText message before sending it to its own children.
This process continues until the leaf agents are reached.

When the leaf agents are reached, the utility phase is initiated.
This second phase is bottom-up and starts from the leaf agents.
Based on the allocated functions, the agents calculate the utility based on the \sampleMessageText message and the assignments of their own variables.
This utility is encoded within a \utilityMessageText message and sent to the parent of the agent.
Upon receiving a \utilityMessageText message, an agent calculates the utility of its utility functions.
The resulting utility value is added to the \utilityMessageText message before sending it to its parent.
This phase finishes when the root agent receives a \utilityMessageText message from all its children.
This moment marks the end of a single iteration within the solver.
At this time, all agents have obtained the utility value associated with the sample of their variables.
This information is used to update the probabilistic measure and thereby the selection of the sample in the next iteration.

\structure{Sample-based solver alternatives}
While the overall procedure is similar for DUCT and \Gibbs, their main difference is in the method of selecting additional samples.
The probabilistic measure in DUCT calculates confidence bounds of the utility of the samples and selects samples to improve these bounds.
The \Gibbs algorithm selects samples based on joint probability distributions and keeps track of the differences between the utility values of the samples as termination criterion.
While \Gibbs is more memory efficient compared to DUCT, both algorithms have a runtime complexity that is linear in the cardinality of the largest domain \cite[Table 4]{fioretto2019distributedreview}.
Therefore, both \Gibbs and DUCT suffer from discretization of continuous domains and are not suitable for continuous DCOPs.

\structure{Consecutive sampling and utility phase result in a problem}
An additional disadvantage of both solvers is the non-determinism with respect to the utility value of a sample.
This is caused by the consecutive sampling and utility phases, since within an iteration all agents sample a single value from their local search space.
In other words, the same \sampleMessageText message can result in different \utilityMessageText messages when the children of an agent select different samples for their variables.

In order to remove the non-determinism, the sampling and utility phase in \DBay will be restricted to parent and children instead of the entire pseudo-tree.
To be more precise, when a child receives a \sampleMessageText message it will first iterate between its own children before sending a \utilityMessageText message to its parent.
This will guarantee that the utility value of the \utilityMessageText message will always be the same for the same \sampleMessageText message.

\structure{Formal algorithm description}
\DBay as described in \cref{algorithm:DBAY} (\cref{section:distributed_bayesian_algorithm}) involves four phases:
\begin{description}[align=left]
    \item[(1) Pseudo-tree construction]
        The agents create a pseudo-tree from the constraint graph of the DCOP, by performing a depth-first search traversal \cite{Awerbuch1985}.
        Afterwards, every agent \agentIndexed knows its parent/children sets (\ParentIndexed/\ChildIndexed) and pseudo-parents/pseudo-children sets (\PseudoParentIndexed/\PseudoChildIndexed), where $\ParentIndexed, \PseudoParentIndexed, \ChildIndexed, \PseudoChildIndexed \subset \agentSet$.
        The pseudo-tree is used as communication structure in which agents only communicate with agents with whom they share a parent/child relation.
        Note that (depending on the DCOP) none of the agents have a complete overview of the pseudo-tree.

    \item[(2) Allocation of utility functions]
        Similar to the allocation of variables, all utility functions in \functionSet are exclusively allocated to the agents.
        Every agent \agentIndexed constructs two separate function sets based on the variables of the agent and the variables of its (pseudo-)parents.
        Firstly, the utility function set \functionSetLocalLongIndexed, which only depends on the agent itself.
        Secondly, the shared utility function set, \functionSetSharedLongIndexed, which involves the agent and its (pseudo-)parents.
        These two function sets are combined as \functionSetLocalCombinedIndexedLong.

    \item[(3) Sample propagation]
        In this phase, every agent optimizes its local variables through the \BO method and the exchange of \sampleMessageText and \utilityMessageText messages.
        By doing so, the assignments of the local variables will converge to the global optimum of the objective function as will be shown in \cref{subsection:globalconvergenceofDBAY}.
        The local variables of \agentIndexed are defined as \localVariablesIndexedLong.
        The variables are optimized based on the aggregate utility of all utility functions in set \functionSetLocalCombinedIndexed and the utility functions of its children ($\utilityFunctionIndexed \in \functionSetLocalCombinedChild$ for all $\agent_\agentIndexChild \in \ChildIndexed$).
        Consequently, the aggregate utility values obtained by the root agent hold the utility values of the objective function.

        Since a sample from (pseudo-)parents is required to calculate the utility of the functions in \functionSetSharedIndexed, every agent \agentIndexed waits for a \sampleMessageText message.
        The phase is therefore initiated by a \sampleMessageText message from the root agent.
        The phase finishes when a convergence threshold is reached by the root agent.
        Upon receiving \sampleMessageText message \sampleMessageIndexedParent from its parent \agentIndexedParent, agent \agentIndexed optimizes its own variables accordingly by sampling all functions in set \functionSetLocalCombinedIndexed.
        The samples are selected through the optimization of an acquisition function.
        Note that the acquisition function is based a kernel and on all preceding samples.

        If the agent is a leaf agent, the agent can optimize its variables without considering the impact of its assignments on other agents.
        However, when the agent has children, the agent needs to send its sample and wait for a \utilityMessageText message to be able to know its effect on the functions in the sets \functionSetLocalCombinedChild  for all $\agent_\agentIndexChild \in \ChildIndexed$.
        Therefore, the agent augments the \sampleMessageText message of its parent with its own sample as \sampleMessageIndexedOtherLong and sends the message to all its children.
        The agent then waits until it has received all \utilityMessageText messages from its children.
        Only then the agent is able to compute the aggregate utility and return a \utilityMessageText message to its own parent.
        Note that the aggregate utility represents the optimal utility for the sample of the agent and all its \mbox{(pseudo-)children}.
        A \utilityMessageText message is defined as \optimumMessageToParentLong, where \mbox{\optimumMessageLocalLong} and \optimumMessageIndexedChildren define the utility and the aggregated child utility, respectively.
        A graphical overview of the \emph{sample phase} is shown in \cref{fig:B-DPOP_sample_phase}.

    \item[(4) Assignment propagation]
        The final phase is the assignment propagation phase, in which the root agent $\agent_\agentIndexRoot$ sends the final assignment of all its variables to its children as a \sampleFinalMessageText message \sampleMessageBestInitialRootLong.
        Based on these assignments the children can assign their own variables to the value corresponding to the optimal utility value.
        Afterwards, every agent adds its assignments to the message as \sampleMessageBestInitialIndexedLong.
        After the leaf agents have received a \sampleFinalMessageText message, all agents have completed their local assignments \assignmentPartialAgentOptimal.
        Note that typically no agent has information of the complete assignment, \mbox{$\assignmentPartialOptimal = \{ \assignmentPartialAgentOptimal ~:~ \variableIndex = 1, \dots, \numberOfAgents \}$}.
\end{description}
\begin{figure}[ht]
    \centering
        \input{figures/Algorithm_sample_phase_example.tikz}
    \caption{
        Graphical overview of the sample phase of \DBay.
        Agents are indicated by circles labeled with an agent index, and utility functions are shown as black lines.
        Starting from the root $\agent_\agentIndexRoot$ (top-left), a \sampleMessageText message $\sampleMessage_\agentIndexRoot$ is sent to its children ($\agent_2, \agent_3$).
        Subsequently, agent $\agent_3$ will send a \sampleMessageText message $\sampleMessage_3$ to its children.
        After iterating between its children and calculating its local utility, agent $\agent_3$ combines all local utilities and sends a \utilityMessageText message $\utilityMessage^\agentIndexRoot_3$ to its parent.
        This process is repeated when the root $\agent_\agentIndexRoot$ sends another  \sampleMessageText message and finishes when $\agent_\agentIndexRoot$ reaches a convergence threshold.
        Note that the interactions between $\agent_\agentIndexRoot$ and $\agent_2$ are not illustrated.
    }
    \label{fig:B-DPOP_sample_phase}
\end{figure}

\clearpage  
\structure{Bridge to analysis section}
In the next section the convergence of \DBay to the global optimum of a continuous DCOP is analysed.
\DBay utilizes \BO for the sample selection within the sample propagation phase.
For that reason, the performance of \DBay depends highly on the properties of the \BO method.
As mentioned in \cref{section:BayesianOptimization}, \BO consists of the combination of a kernel and an acquisition function.
Therefore, the analysis is focussed on the selection of the kernel, the acquisition function, and their parameters.

\section{Theoretical analysis of \DBay}
\label{section:analysis}

\structure{General introduction}
This section will analyse the convergence of \DBay to the global optimum of a continuous DCOP in two parts.
Firstly, in \cref{subsection:convergenceofBO}, the convergence to the global optimum of the utility functions within the sampling phase is proven.
It will be shown that if the Lipschitz constant of the utility functions is known, the convergence to the global optimum can be guaranteed through the appropriate selection of a kernel and an acquisition function.
In this work, all utility functions are assumed to be Lipschitz continuous with known Lipschitz constant.
A utility function \utilityFunctionDot is Lipschitz continuous with Lipschitz constant \lipschitzConstantIndexed if
\begin{align}
    \lipschitzConstantIndexedLong \qquad \forall \functionInputBaseIndexed, \functionInputBaseIndexedOther \in \domainOperator(\utilityFunction)
    \label{eq:LipschitzConstant}
\end{align}
where $\domainOperator(\utilityFunction)$ denotes the domain of the utility function.

Secondly, in \cref{subsection:globalconvergenceofDBAY}, the convergence of \DBay to the global optimum of the objective function based on the global optima of the utility functions is proven.
This analysis focusses on the assignment propagation phase.
The two parts of the analysis are combined to prove convergence of \DBay to the global optimum of continuous DCOPs with utility functions with known Lipschitz constants.

\subsection{Convergence of \BO based on Lipschitz continuous functions}
\label{subsection:convergenceofBO}

\structure{Local optimization through dense sampling of the domains}
As shown by \citeA{torn1989global}, the convergence to the global optimum of a function by \BO can only be guaranteed through \emph{dense} sampling of the domain of the function.
For this reason, within the \BO method, the acquisition function will need to produce dense samples.
In the work of \citeA[Theorem~6]{Vazquez2010}, the \EI acquisition function, given by \cref{eq:acquisition}, is proven to produce dense observations within its search region.
The search region is defined in Definition \ref{def:searchRegion}.
\begin{defn}
    \label{def:searchRegion}
    The search region of the expected improvement acquisition function $\acquisitionFunctionBase(\cdot, \acquisitionFunctionParameter | \observationSet)$ (based on \utilityFunctionDot and \observationSet) is defined by
    \begin{align*}
        \searchRegionLong.
    \end{align*}
\end{defn}

\structure{how to achieve dense sampling? -> how to ensure $S$ includes $x^*$?}
As a consequence of the dense sample generation property, \searchRegion will converge to an empty set when the number of samples goes to infinity.
Therefore, since the next sample is chosen from the search region ($\nextSample \in \searchRegion$), the global optimum will be sampled for $\observationIndex \rightarrow \infty$ if $\functionInputOptimal \in \searchRegion$.
In other words, if the optimum inclusion ($\functionInputOptimal \in \searchRegion$) property holds, then global convergence is guaranteed.

\begin{figure}[ht]
    \centering
    \includegraphics[width=0.85\textwidth, trim={0px 0px 0px 0px}, clip]{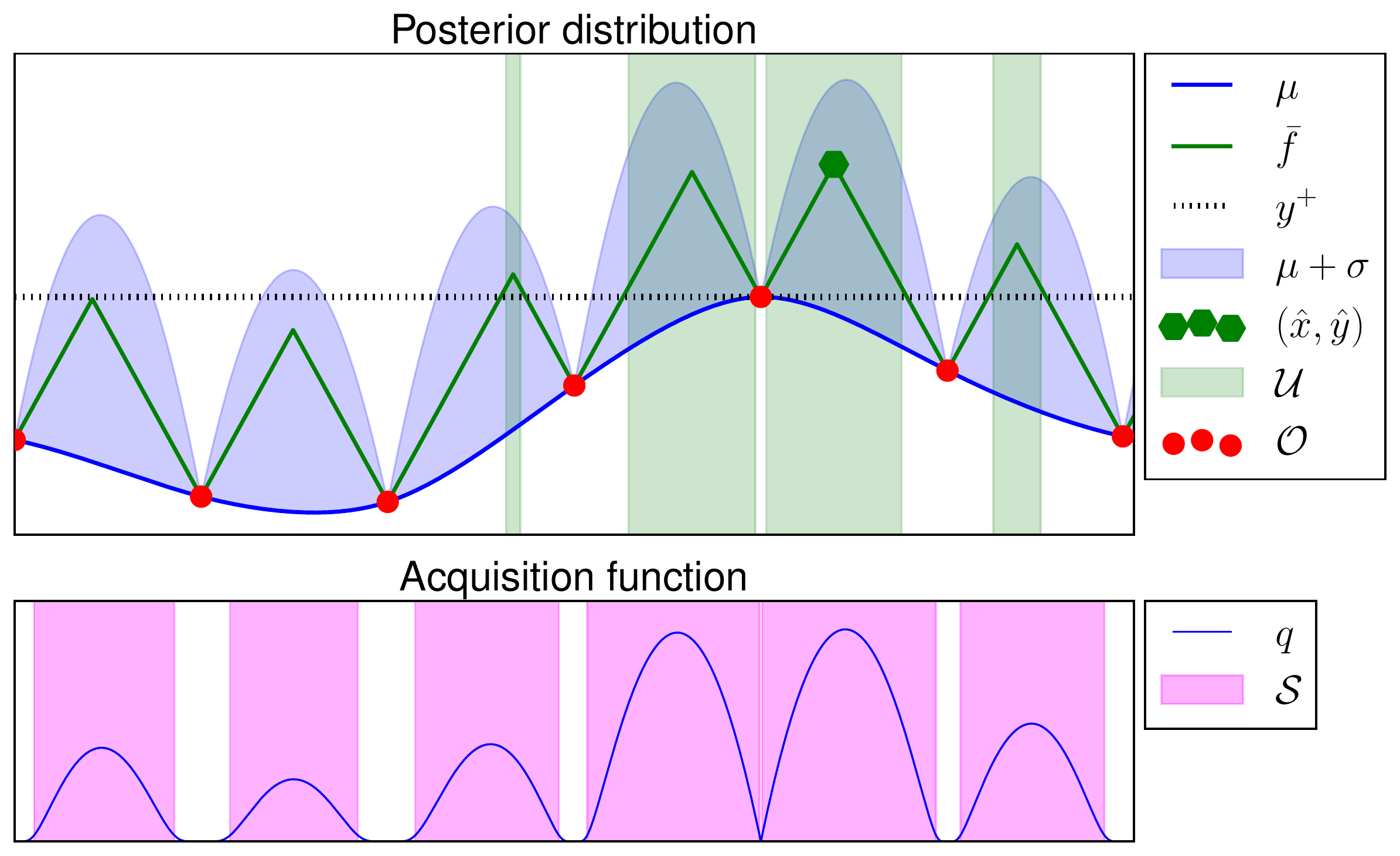}
    \caption{
        Graphical overview of the sets \upperBoundRegion, \searchRegion, and the upper bound function \upperBoundFunctionDot based on the observations \observationSet.
        Based on the observations (red circles), the corresponding upper bound function \upperBoundFunctionBase, the posterior mean function \meanFunctionBase and the the upper bound region \upperBoundRegion are shown.
        }
        \label{fig:LBEI}
\end{figure}

\structure{Based on upper bound function the set $U$ includes $x^*$}
In order to find the conditions for which the optimum inclusion holds, the upper bound region is introduced.
The upper bound region set \upperBoundRegion (Definition \ref{def:upperBoundRegion}) is based on the upper bound function \upperBoundFunctionDot (Definition \ref{def:upper_bound_function}).
Note that by definition, $\upperBoundFunctionBase(\functionInputBase) \geq \utilityFunction(\functionInputBase)$ for all $\functionInputBase$.
As shown in Lemma \ref{lemma:optimumInU}, this region is guaranteed to include the global optimum if the optimum has not already been observed.
A graphical example of \upperBoundRegion, \searchRegion, and \upperBoundFunctionDot can be seen in \cref{fig:LBEI}.

\begin{defn}
    \label{def:upper_bound_function}
    The upper bound function \upperBoundFunction of (a Lipschitz continuous) function \utilityFunctionDot over all observations in \observationSet is defined by
    \begin{align*}
        \upperBoundFunctionLongIndex \qquad \forall \functionInputBase \in \domainOperator(\utilityFunction).
    \end{align*}
\end{defn}
\begin{defn}
    \label{def:upperBoundRegion}
    The upper bound region of \utilityFunctionDot holds all values of \functionInputBase for which the upper bound function \upperBoundFunctionDot (Definition \ref{def:upper_bound_function}) is larger than the maximum observed value \observationMaxSample and is defined by
    \begin{align*}
        \upperBoundRegionLong.
    \end{align*}
\end{defn}
Note that by definition \upperBoundRegion does not include any observations in \observationSet since $\upperBoundFunction = \observationOutputIndexed \leq \observationMaxSample$ for all \observationInputIndexed.

\begin{lemma}
    \label{lemma:optimumInU}
    The upper bound region \upperBoundRegion includes the optimum sample \functionInputOptimal if it has not been observed.
    Formally, if $\observationMaxSampleInput \neq \functionInputOptimal$, then $\functionInputOptimal \in \upperBoundRegion$.
    \begin{proof}
        By Definition \ref{def:upper_bound_function} and \cref{eq:LipschitzConstant}, the value of the upper bound function at the optimal sample is larger than the optimal value.
        Formally, $\bar{\utilityFunction}(\functionInputOptimal | \observationSet) \geq \functionOutputOptimal$.
        If the optimal sample has not been observed ($\functionInputOptimal \neq \observationMaxSampleInput$), then $\functionOutputOptimal > \observationMaxSample$.
        Consequently, $\bar{\utilityFunction}(\functionInputOptimal | \observationSet) \geq \observationMaxSample$.
        Therefore, by definition of the upper bound region (Definition \ref{def:upperBoundRegion}) the optimal sample is included ($\functionInputOptimal \in \upperBoundRegion$).
    \end{proof}
\end{lemma}

\structure{how can we relate $U$ and $S$?}
Based on the definition of the upper bound region, the optimum inclusion is satisfied when the set inclusion $\upperBoundRegion \subseteq \searchRegion$ holds.
The conditions for the set inclusion are given in two parts.
Firstly, in Lemma \ref{lemma:x_in_searchRegion} it is shown that if $\meanFunction + \deviationFunction \geq \observationMaxSample + \acquisitionFunctionParameter$ and $\deviationFunction > 0$ then $\functionInputBase \in \searchRegion$.
Secondly, by Definition \ref{def:upperBoundRegion}, if $\upperBoundFunction > \observationMaxSample$ then $\functionInputBase \in \upperBoundRegion$.
By combining both conditions, we find $\upperBoundRegion \subseteq \searchRegion$ if $\meanFunction + \deviationFunction \geq \upperBoundFunction + \acquisitionFunctionParameter$ for all $\functionInputBase \in \domainOperator(\utilityFunction)$, as shown in Lemma \ref{lemma:setInclusion}.

\newcommand{\w}{
    \mathcmd{
        w(\functionInputBase, \acquisitionFunctionParameter | \observationSet)
    }
}
\newcommand{\z}{
    \mathcmd{
        z(\functionInputBase, \acquisitionFunctionParameter | \observationSet)
    }
}
\newcommand{\h}{
    \mathcmd{
        h(z)
    }
}
\begin{lemma}
    \label{lemma:x_in_searchRegion}
    If $\meanFunction + \deviationFunction \geq \observationMaxSample + \acquisitionFunctionParameter$ and $\deviationFunction > 0$ then $\functionInputBase \in \searchRegion$.
    \begin{proof}
        By Definition \ref{def:searchRegion}, $\functionInputBase \in \searchRegion$ if $\acquisitionFunction > 0$.
        Define $\w = \acquisitionHelperExtended$ and $\z = \frac{\w}{\deviationFunction}$ and through substitution rewrite \cref{eq:acquisition} as
        \begin{align}
            \acquisitionFunction &= \w \Phi\bracketBig{\z} + \deviationFunction \phi\bracketBig{\z} \text{ when } \deviationFunction > 0.
            \label{lemma:eq:acquisition_function}
        \end{align}
        Note that $\acquisitionFunction = 0$ when $\deviationFunction = 0$.
        \\
        \\
        Let $\meanFunction + \deviationFunction \geq \observationMaxSample + \acquisitionFunctionParameter$, then $\w \geq -\deviationFunction$.
        \\
        \\
        To show that $\acquisitionFunction > 0$ when $\w \geq -\deviationFunction$, two separate cases for a given $\functionInputBase \in \domainOperator(\utilityFunction)$ are considered: $\w \geq 0$ and $\w < 0$:
        \begin{description}
            \item[Case 1:]
            Let $\w \geq 0$. Then $\z \geq 0$ since $\deviationFunction > 0$ for all $\functionInputBase \in \domainOperator(\utilityFunction)$.
            Since $\Phi\bracket{z} > 0$ and $\phi\bracket{z} > 0$ for $z \geq 0$, through substitution in \cref{lemma:eq:acquisition_function} we find $\acquisitionFunction > 0$.

            \item[Case 2:]
            Let $\w < 0$, Then $-\w < \deviationFunction$ and subsequently $\z$ is bounded to the interval $[-1, 0]$.
            This interval is analysed by applying \cref{lemma:eq:acquisition_function} as
            \begin{align*}
                \acquisitionFunction &> 0
                \\
                \w \Phi\bracketBig{\z} + \deviationFunction \phi\bracketBig{\z} &> 0
                \\
                \z \Phi\bracketBig{\z} + \phi\bracketBig{\z} &> 0
            \end{align*}
            Define $\h = z \Phi\bracket{z} + \phi\bracket{z}$. Then since $\Phi^{\prime}\bracket{z} = \phi\bracket{z}$ and $\phi\bracket{z} = \frac{1}{\sqrt{2 \pi}} \mathrm{e}^{-\frac{1}{2} z^2}$, we find $h^{\prime}(z) = \Phi\bracket{z} + z \phi\bracket{z} - z \phi\bracket{z} = \Phi\bracket{z}$.
            For $z$ in the interval $[-1, 0]$ we find
            \begin{align*}
                \h
                &= \int_{-1}^{z} h^{\prime}(v) \mathrm{d}v + h(-1)
                = \int_{-1}^{z} \Phi(v) \mathrm{d}v -\Phi(-1) + \phi(-1) > 0
            \end{align*}
            since $\Phi(z) > 0$ for finite inputs, and $-\Phi(-1) + \phi(-1) > 0$.
            Therefore, we conclude that $\acquisitionFunction > 0$.
        \end{description}
        In both cases we find $\acquisitionFunction > 0$.
        Therefore, if $\meanFunction + \deviationFunction \geq \observationMaxSample + \acquisitionFunctionParameter$ and $\deviationFunction > 0$, then $\functionInputBase \in \searchRegion$.
    \end{proof}
\end{lemma}

\begin{lemma}
    If $\meanFunction + \deviationFunction \geq \upperBoundFunction + \acquisitionFunctionParameter$ for all $\functionInputBase \in \domainOperator(\utilityFunction)$, then $\upperBoundRegion \subseteq \searchRegion$.
    \begin{proof}
        As shown in Lemma \ref{lemma:x_in_searchRegion}, if $\meanFunction + \deviationFunction \geq \observationMaxSample + \acquisitionFunctionParameter$ and $\deviationFunction > 0$, then $\functionInputBase \in \searchRegion$.
        By definition of \upperBoundRegion, for all $\functionInputBase \in \upperBoundRegion$ we find $\upperBoundFunction > \observationMaxSample$.
        Additionally, for all $\functionInputBase \in \upperBoundRegion$ we find $\deviationFunction > 0$, since $\varianceFunction = 0$ only if $\functionInputBase = \observationFunctionInput$ and $\observationFunctionInput \notin \upperBoundRegion$.
        Therefore, if $\meanFunction + \deviationFunction \geq \upperBoundFunction + \acquisitionFunctionParameter$, then $\functionInputBase \in \searchRegion$ for all $\functionInputBase \in \upperBoundRegion$.
        In conclusion, $\upperBoundRegion \subseteq \searchRegion$.
    \end{proof}
\label{lemma:setInclusion}
\end{lemma}

\structure{how to ensure the inequality holds?}
Note that \meanFunctionDot and \deviationFunctionDot depend on the kernel and \upperBoundFunctionDot depends on the Lipschitz constant.
This raises the question of which (type of) kernel is capable of explicitly associating its mean function and variance function to the Lipschitz constant.

\structure{Use Markovian kernel for its simplicity/similarity to upper bound function}
An answer can be found in the work of \citeA{Ding2018}.
\citeA{Ding2018} introduced the Markovian class kernels as kernels that possess a Markovian posterior distribution.
A Markovian posterior distribution for a certain input only depends on the observations surrounding that input.
This property is beneficial as the upper bound function, which is directly related to the Lipschitz constant, possesses the same property.
An additional benefit of this class of kernels is that the elements of $\kernelGramian^{-1}(\observationSet)$ can be expressed analytically.
This removes the need of inversion of a matrix of which the size grows with the number of observations, since $\kernelGramian(\observationSet) \in \realsSet^{\numberOfObservations \times \numberOfObservations}$.
As noted by \citeA{Rasmussen2004}, this inversion is considered a major restriction to the practical application of \BO.
\noindent
In general, a Markovian class kernel is defined by
\begin{align*}
    \kernel(\functionInputBaseIndexed, \functionInputBaseIndexedOther) = \kernelScale^2 \Big( p(\functionInputBaseIndexed)g(\functionInputBaseIndexedOther) \mathds{I}_{\functionInputBaseIndexed \leq \functionInputBaseIndexedOther} +p(\functionInputBaseIndexedOther)g(\functionInputBaseIndexed) \mathds{I}_{\functionInputBaseIndexed > \functionInputBaseIndexedOther} \Big)
\end{align*}
for some function $p(\cdot)$ and $g(\cdot)$, where $\mathds{I}(\cdot)$ is the indicator function and \kernelScale is the kernel scale parameter.
The mean function \meanFunctionObservationDot and the variance function $\varianceFunctionBase_\observationIndex(\cdot | \observationSet)$ of the posterior on the interval between observations for a kernel of this class is defined by,
\begin{align}
    \meanFunctionObservation &=
        \kernelBold^\transpose_\observationIndex (\functionInputBase, \observationSet)
        \kernelGramian^{-1}_\observationIndex (\observationSet)
        \functionOutput_\observationIndex (\observationSet)
    \label{eq:mean_function_interval}
    \\
    \varianceFunctionBase_\observationIndex(\functionInputBase | \observationSet) &=
    \kernel\left( \functionInputBase, \functionInputBase \right) - \kernelBold_\observationIndex^\transpose (\functionInputBase, \observationSet)
    \kernelGramian^{-1}_\observationIndex(\observationSet)
    \kernelBold_\observationIndex (\functionInputBase, \observationSet)
    \label{eq:variance_function_interval}
\end{align}
for $\functionInputBase \in [\functionInputObservationLeft,  \functionInputObservation]$, where
\begin{align*}
    \kernelBold_\observationIndex(\functionInputBase, \observationSet)&=
    \begin{bmatrix}
        \kernel\bracket{\functionInputBase_1, \functionInputBase} &
        \dots &
        \kernel\bracket{\functionInputObservationLeft, \functionInputBase} &
        \kernel\bracket{\functionInputObservation, \functionInputBase} &
        \kernel\bracket{\functionInputObservationRight, \functionInputBase} &
        \dots &
        \kernel\bracket{\functionInputBase_\numberOfObservations, \functionInputBase}
    \end{bmatrix}^\transpose,
    \\
    \functionOutput_\observationIndex(\observationSet) &=
    \begin{bmatrix}
        \functionOutputBase_1 &
        \dots &
        \functionOutputObservationLeft &
        \functionOutputObservation &
        \functionOutputObservationRight &
        \dots &
        \functionOutputBase_\numberOfObservations
    \end{bmatrix}^\transpose,
\end{align*}
and $\kernelGramian^{-1}_\observationIndex(\observationSet)$ is a tridiagonal matrix of appropriate dimensions where the tridiagonal elements of $\kernelGramian_\observationIndex^{-1}(\observationSet)$, for $\numberOfObservations \geq 3$ and if $\kernelGramian_\observationIndex(\observationSet)$ is nonsingular, are given by
\begin{align*}
    (\kernelGramian_\observationIndex^{-1}(\observationSet))_{\observationIndex,\observationIndex}
    =
    \begin{cases}
        \frac{
            \kernelScale^{-2}
            p(\functionInputBase_2)
        }{
            p(\functionInputBase_1)
            \bracketbig{
                p(\functionInputBase_2) g(\functionInputBase_1) - p(\functionInputBase_1) g(\functionInputBase_2)
            }
        }, &\text{if } \observationIndex = 1,
        \\[12pt]
        \frac{
            \kernelScale^{-2}
            \bracketbig{
                p(\functionInputBase_{\observationIndex+1}) g(\functionInputBase_{\observationIndex-1}) - p(\functionInputBase_{\observationIndex-1})g(\functionInputBase_{\observationIndex+1})
            }
        }{
            \bracketbig{
                p(\functionInputBase_{\observationIndex}) g(\functionInputBase_{\observationIndex-1}) - p(\functionInputBase_{\observationIndex-1}) g(\functionInputBase_{\observationIndex})
            }
            \bracketbig{
                p(\functionInputBase_{\observationIndex+1})  g(\functionInputBase_{\observationIndex}) - p(\functionInputBase_{\observationIndex}) g(\functionInputBase_{\observationIndex+1})
            }
        }, &\text{if } \observationIndex \in \{ 2, \dots, \numberOfObservations - 1\},
        \\[12pt]
        \frac{
            \kernelScale^{-2}
            g(\functionInputBase_{\numberOfObservations-1})
        }{
            g(\functionInputBase_{\numberOfObservations})
            \bracketbig{
                p(\functionInputBase_{\numberOfObservations}) g(\functionInputBase_{\numberOfObservations-1}) - p(\functionInputBase_{\numberOfObservations - 1}) g(\functionInputBase_{\numberOfObservations})
            }
        }, &\text{if } \observationIndex = \numberOfObservations,
    \end{cases}
\end{align*}
and
\begin{align*}
    (\kernelGramian_\observationIndex^{-1}(\observationSet))_{\observationIndex-1,\observationIndex}
    =
    (\kernelGramian_\observationIndex^{-1}(\observationSet))_{\observationIndex,\observationIndex-1}
    =
    \frac{
        -\kernelScale^{-2}
    }{
        \bracketbig{
            p(\functionInputBase_{\observationIndex}) g(\functionInputBase_{\observationIndex-1}) - p(\functionInputBase_{\observationIndex-1}) g(\functionInputBase_{\observationIndex})
        }
    }
    ,~\observationIndex = 2,\dots,\numberOfObservations.
\end{align*}
All other elements of $\kernelGramian_\observationIndex^{-1}(\observationSet)$ are equal to zero.

\structure{which kernel and what parameters to choose?}
Next, we show that for the Dirichlet kernel, as introduced by \citeA{Ding2018}, the inequality of Lemma \ref{lemma:setInclusion} will hold for all observations if the kernel scale is chosen appropriately.
This kernel is selected over other Markovian class kernels because of to its simplicity.
The Dirichlet kernel defined by
\begin{align}
    \label{eq:dirichlet}
    \dirichletKernelLong
\end{align}
for $\functionInputBaseIndexed, \functionInputBaseIndexedOther \in [0, 1]$.
Note that for \dirichletKernel, we find that $p(\functionInputBase) = \functionInputBase$ and $g(\functionInputBase) = (1 - \functionInputBase)$.
The mean function, given by \cref{eq:mean_function_interval}, and the variance function, given by \cref{eq:variance_function_interval}, corresponding to the Dirichlet kernel in the interval $[\functionInputObservationLeft, \functionInputObservation]$ can be written as
\begin{align}
    \meanFunctionObservation &=
    \frac{\functionOutputBase_{\observationIndex-1}(\functionInputObservation - x) + \functionOutputBase_{\observationIndex}(x - \functionInputObservationLeft)}{\functionInputObservation - \functionInputObservationLeft},
    \label{eq:mean_function_interval_Dirichlet}
\end{align}
\begin{align}
    \varianceFunctionBase_\observationIndex(\functionInputBase | \observationSet)
    &=
        \kernelScale^2
        \frac{
            -(\functionInputObservation - \functionInputBase)
            (\functionInputObservationLeft - \functionInputBase)
        }{
            \functionInputObservation - \functionInputObservationLeft
        }.
    \label{eq:variance_function_interval_Dirichlet}
\end{align}
The derivation of \cref{eq:mean_function_interval_Dirichlet,eq:variance_function_interval_Dirichlet} can be found in \cref{section:dirichlet_interval_functions}.
Note that both the mean function and the variance function on the interval $[\functionInputObservationLeft, \functionInputObservation]$ only depend on the observations ($\observation_{\observationIndex-1}$, $\observation_{\observationIndex}$) at the boundaries of the interval.

\structure{Sufficient to proof for a single segment}
Based on \cref{eq:mean_function_interval_Dirichlet,eq:variance_function_interval_Dirichlet}, the inequality of Lemma \ref{lemma:setInclusion} holds if
\begin{align}
    \meanFunctionObservation  + \deviationFunctionObservation
    &\geq
    \upperBoundFunction + \acquisitionFunctionParameter
    \text{ for }\functionInputBase \in [\functionInputObservationLeft, \functionInputObservation]
    \label{eq:inequality_observation}
\end{align}
for all $\observationIndex \in \{1, \dots, \numberOfObservations\}$, given $\functionInputBase_1 = 0$ and $\functionInputBase_\numberOfObservations = 1$.
In other words, by using the Dirichlet kernel, instead of analyzing the inequality in Lemma \ref{lemma:setInclusion} over the entire domain of the function, it is sufficient to analyze \cref{eq:inequality_observation} on the intervals between the observations.

\structure{Theorem (to conclude the section)}
Now that the acquisition function and the kernel have been selected, we need to find their parameters (\acquisitionFunctionParameter and \kernelScale) such that the inequality of \cref{eq:inequality_observation} holds.
These parameters can be appropriately chosen based on the Lipschitz constant as follows:

\begin{theorem}
\label{theorem:kernel_scale}
    For a function $\utilityFunction(\cdot)$ with known Lipschitz constant $\lipschitzConstantIndexed$ and $\domainOperator(\utilityFunction) = [0, 1]$, the Dirichlet kernel \dirichletKernel, and $\numberOfObservations \geq 3$, where $\functionInputBase_1 = 0$ and $\functionInputBase_\numberOfObservations = 1$,
    will yield $\meanFunction + \deviationFunction \geq \upperBoundFunction$ for all $\functionInputBase \in \domainOperator(\utilityFunction)$ if $\kernelScale \geq \lipschitzConstantIndexed$.
    \begin{proof}
        Let the functions \meanFunctionObservationDot and \deviationFunctionObservationDot be as defined by \cref{eq:mean_function_interval_Dirichlet,eq:variance_function_interval_Dirichlet}, respectively.
        At the observations ($\functionInputBase = \functionInputObservation$ for $\observationIndex \in \{1, \dots, \numberOfObservations\}$), the inequality $\meanFunction + \deviationFunction \geq \upperBoundFunction$ is satisfied, since $\meanFunctionBase_\observationIndex(\functionInputObservation | \observationSet)  + \deviationFunctionBase_\observationIndex(\functionInputObservation | \observationSet) = \upperBoundFunctionBase(\functionInputObservation | \observationSet) = \functionOutputObservation$.
        Therefore, by letting $\functionInputBase_1 = 0$ and $\functionInputBase_\numberOfObservations = 1$, only the closed intervals $\functionInputBase \in [\functionInputObservationLeft, \functionInputObservation]$ for all $\observationIndex \in \{2, \dots, \numberOfObservations\}$ need to be examined.
        The proof will focus on these closed intervals next.

        \structure{domain to interval}
        Based on \cref{eq:mean_function_interval_Dirichlet,eq:variance_function_interval_Dirichlet} the inequality $\meanFunction + \deviationFunction \geq \upperBoundFunction$ for the Dirichlet kernel at the closed intervals can be rewritten as
        \begin{align}
            \meanFunctionObservation  + \deviationFunctionObservation
            &\geq
            \upperBoundFunction
            \text{ for }\functionInputBase \in [\functionInputObservationLeft, \functionInputObservation]
            \label{eq:inequality_observation_theorem}
        \end{align}
        \structure{Not at the observations}
        For the benefit of the analysis, we normalize the function input for every interval by defining a normalized function argument as
        \begin{align}
            \functionInputObservationNormalized &=
            \frac{
                \functionInputBase - \functionInputObservationLeft
            }{
                \functionInputObservation - \functionInputObservationLeft
            }.
            \label{eq:input_normalized}
        \end{align}
        Then $\functionInputObservationNormalized \in [0,1]$ for $\functionInputBase \in [\functionInputObservationLeft, \functionInputObservation]$.
        Additionally, define $\functionInputObservationDelta = \functionInputObservation - \functionInputObservationLeft$ as the size of the interval.
        Likewise, the upper bound function \upperBoundFunctionDot can be rewritten based on the normalized interval as
        \begin{align}
            \upperBoundFunctionObservation &=
            \begin{cases}
                \lipschitzConstantNormalized \functionInputObservationNormalized + \functionOutputObservationLeft
                    &\text{if }
                    0  \leq \functionInputObservationNormalized \leq \hat{\functionInputObservationNormalizedBase}_{\observationIndex} \\
                \lipschitzConstantNormalized (1 - \functionInputObservationNormalized) + \functionOutputObservation
                    &\text{if }
                    \hat{\functionInputObservationNormalizedBase}_{\observationIndex} < \functionInputObservationNormalized \leq 1
            \end{cases}
            \label{eq:upper_bound_observation}
        \end{align}
    By defining an auxillary variable $\functionOutputObservationConstant \in [-1, 1]$ as
    \begin{align*}
        \functionOutputObservationConstant
        &=
        \frac{\functionOutputObservation - \functionOutputObservationLeft}{\lipschitzConstantIndexed (\functionInputObservation - \functionInputObservationLeft)}
        =
        \frac{\functionOutputObservation - \functionOutputObservationLeft}{\lipschitzConstantNormalized}
    \end{align*}
    all possible critical points of \upperBoundFunctionDot ($\hat{\functionInputObservationNormalizedBase}_{\observationIndex}$) can be written as a function of $\functionOutputObservationConstant$ as
        \begin{align}
            \functionInputObservationNormalizedCriticalPoint &=
            \frac{1}{2}
            +
            \frac{\bracket{\functionOutputObservation - \functionOutputObservationLeft}}{2 \lipschitzConstantNormalized}
            \\&=
            \frac{1}{2}
            \left(
                1 + \functionOutputObservationConstant
            \right)
            .
            \label{eq:CriticalPoint}
        \end{align}
        A graphical overview of the normalized interval and the corresponding functions can be seen in \cref{fig:overview_proof}.
        \begin{figure}[ht]
            \centering
            \includegraphics[width=0.60\textwidth, trim={0px 0px 0px 0px}, clip]{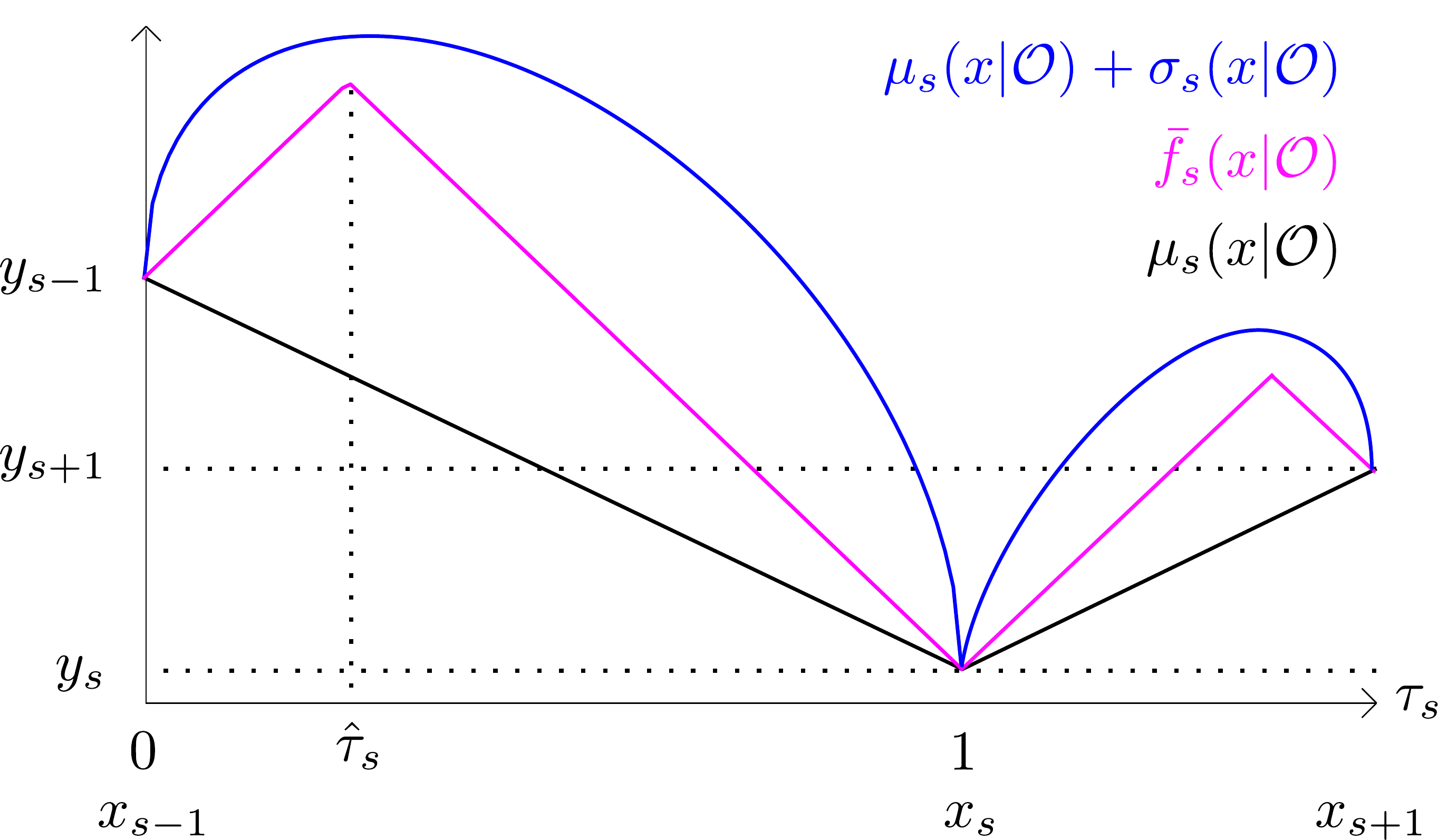}
            \caption{
                Graphical overview of the interval $[\functionInputObservationLeft, \functionInputObservationRight]$ where \meanFunctionObservationDot and \deviationFunctionObservationDot are based on the Dirichlet kernel and \upperBoundFunctionObservationDot is based on \lipschitzConstantNormalized.
                The normalized function argument $\functionInputObservationNormalized$ is shown over its domain $[0, 1]$.
                Based on the upper bound function, the domain can be divided into two intervals $[0, \hat{\functionInputObservationNormalizedBase}_{\observationIndex})$ and $[\hat{\functionInputObservationNormalizedBase}_{\observationIndex}, 1]$.
                \label{fig:overview_proof}
            }
        \end{figure}

        \structure{After defining these variables, we can write the mean and deviation for the normalized interval}
        After defining these variables, two separate intervals can be considered, $[0, \functionInputObservationNormalizedCriticalPoint)$ and $[\functionInputObservationNormalizedCriticalPoint, 1]$.
        Both these intervals will be investigated next.
        \begin{description}
            \item[Interval $\ensuremath{[0, \functionInputObservationNormalizedCriticalPoint)}$:] 
                Let $\functionInputObservationNormalized \in [0, \functionInputObservationNormalizedCriticalPoint)$.
                The mean function, given by \cref{eq:mean_function_interval_Dirichlet}, on the interval can be rewritten based on the normalized function argument as
                \begin{align}
                    \meanFunctionObservation
                    &=
                    \frac{\functionOutputBase_{\observationIndex-1}(\functionInputObservation - x) + \functionOutputBase_{\observationIndex}(x - \functionInputObservationLeft)}{\functionInputObservation - \functionInputObservationLeft}
                    \nonumber \\
                    &=
                    \functionOutputBase_{\observationIndex-1}(1 - \functionInputObservationNormalized) +
                    \functionOutputBase_{\observationIndex} \functionInputObservationNormalized
                    \nonumber \\
                    &=
                    \functionOutputBase_{\observationIndex-1}
                    +
                    \functionInputObservationNormalized
                    \lipschitzConstantNormalized
                    \functionOutputObservationConstant.
                    \label{eq:mean_function_interval_left}
                \end{align}

                Likewise the varince function in \cref{eq:variance_function_interval_Dirichlet}) can be rewritten as
                \begin{align}
                    \deviationFunctionObservation
                    &=
                    \kernelScale
                    \sqrt{
                        \frac{
                            -(\functionInputObservation - \functionInputBase)
                            (\functionInputObservationLeft - \functionInputBase)
                        }{
                            \functionInputObservation - \functionInputObservationLeft
                        }
                    }
                    \nonumber \\
                    &=
                    \kernelScale
                    \sqrt{
                        (1 - \functionInputObservationNormalized)
                        \functionInputObservationNormalized
                        \Delta \functionInputBase_\observationIndex
                    }.
                    \label{eq:deviation_function_interval_left}
                \end{align}

                Substitution of \cref{eq:mean_function_interval_left,eq:deviation_function_interval_left,eq:upper_bound_observation} into \cref{eq:inequality_observation_theorem} yields
                \begin{align}
                    \meanFunctionObservation  + \deviationFunctionObservation &\geq \upperBoundFunctionObservation
                    \nonumber
                    \\
                    \kernelScale
                    &\geq
                    (1 - \functionOutputObservationConstant)
                    \sqrt{
                        \frac{
                            \functionInputObservationNormalized
                        }{
                            (1 - \functionInputObservationNormalized)
                        }
                    }
                    \lipschitzConstantIndexed
                    \sqrt{\Delta \functionInputBase_\observationIndex}
                    .
                    \label{eq:theorem_interval_left}
                \end{align}
                Note that \cref{eq:theorem_interval_left} gives an explicit expression of the value for the kernel scale $\kernelScale$ and all possible values of input/output pairs of the observations through the auxillary variable $\functionOutputObservationConstant$.

                \structure{The next step is to find the bounds for the parts of the RHS}
                To analyze the values of $\kernelScale$ for which the inequality of \cref{eq:theorem_interval_left} holds, the upper bound of the right-hand side is determined.

                Since $\sqrt{\functionInputObservationNormalized / (1 - \functionInputObservationNormalized)}$ is monotonically increasing with respect to \functionInputObservationNormalized, we find for $\functionInputObservationNormalized$ in the interval $[0, \functionInputObservationNormalizedCriticalPoint)$,
                \begin{align*}
                (1 - \functionOutputObservationConstant)
                \sqrt{\frac{\functionInputObservationNormalized}{(1 - \functionInputObservationNormalized)}}
                &\leq
                (1 - \functionOutputObservationConstant)
                \sqrt{\frac{\functionInputObservationNormalizedCriticalPoint}{(1 - \functionInputObservationNormalizedCriticalPoint)}}
                \\&\leq
                (1 - \functionOutputObservationConstant)
                \sqrt{
                    \frac{
                        \frac{1}{2} (1 + \functionOutputObservationConstant)
                    }{
                        (1 - \frac{1}{2} (1 + \functionOutputObservationConstant))
                    }
                }
                \\&\leq
                    \sqrt{ 1 - \functionOutputObservationConstant^2 }
                \\&\leq 1.
                \end{align*}
                Furthermore, since $\Delta \observationFunctionInput \leq \Delta \functionInputBase \leq 1$, we find through substitution of the upper bounds in \cref{eq:theorem_interval_left} that if $\kernelScale \geq \lipschitzConstantIndexed$, then \cref{eq:theorem_interval_left} is satisfied for all possible observations.

            \item[Interval $\ensuremath{[\functionInputObservationNormalizedCriticalPoint, 1]}$:]
                For this interval the same approach is applied.
                Let $\functionInputObservationNormalized \in [\functionInputObservationNormalizedCriticalPoint, 1]$; then substitution of \cref{eq:mean_function_interval_left,eq:deviation_function_interval_left,eq:upper_bound_observation} in \cref{eq:inequality_observation_theorem} yields
                \begin{align}
                    \meanFunctionObservation  + \deviationFunctionObservation &\geq \upperBoundFunctionObservation
                    \nonumber
                    \\
                    \kernelScale
                    &\geq
                    (1 + \functionOutputObservationConstant)
                    \sqrt{
                        \frac{(1 - \functionInputObservationNormalized)}{\functionInputObservationNormalized}
                    }
                    \lipschitzConstantIndexed \sqrt{\Delta \observationInputIndexed}
                    .
                    \label{eq:theorem_interval_right}
                \end{align}

                Since $\sqrt{(1 - \functionInputObservationNormalized)/ \functionInputObservationNormalized}$ is monotonically decreasing with respect to \functionInputObservationNormalized, we find
                \begin{align}
                    (1 + \functionOutputObservationConstant)
                    \sqrt{
                        \frac{(1 - \functionInputObservationNormalized)}{\functionInputObservationNormalized}
                    }
                    &\leq
                    1.
                \end{align}
                Therefore, we conclude that for the interval $[\functionInputObservationNormalizedCriticalPoint, 1]$ if $\kernelScale \geq \lipschitzConstantIndexed$ then \cref{eq:theorem_interval_right} is satisfied for all possible observations.

            \end{description}
        In conclusion, if $\kernelScale \geq \lipschitzConstantIndexed$, we find that the inequality of \cref{eq:inequality_observation_theorem} will hold for the intervals $[\functionInputObservationLeft, \functionInputObservation]$ for $\observationIndex \in \{1, \dots, \numberOfObservations\}$.
        Since $\functionInputBase_1 = 0$ and $\functionInputBase_\numberOfObservations = 1$, \cref{eq:inequality_observation} hold for all $\functionInputBase \in (0, 1)$.
        Subsequently, $\meanFunction + \deviationFunction \geq \upperBoundFunction$ will holds for all $\functionInputBase \in [0, 1]$.
    \end{proof}
\end{theorem}

According to Theorem \ref{theorem:kernel_scale}, if $\kernelScale \geq \lipschitzConstantIndexed$ we find $\meanFunction + \deviationFunction \geq \upperBoundFunction$ for all \mbox{$\functionInputBase \in \domainOperator(\utilityFunction)$}.
Applying Lemma \ref{lemma:setInclusion} and setting $\acquisitionFunctionParameter = 0$, yields $\upperBoundRegion \subseteq \searchRegion$ where $\functionInputOptimal \in \searchRegion$ for all observations.
Subsequently, the Bayesian optimization will converge to the global optimum of the function.
Note that for all functions without normalized domain, the Lipschitz constant should be scaled according to the scaling required for normalization of the domain.

\subsection{Convergence of \DBay based on global optima of utility functions}
\label{subsection:globalconvergenceofDBAY}
\structure{We go by the assumption that the agents are able to find the optimal assignments}
As shown in \cref{subsection:convergenceofBO}, all agents are able to find the global optimum of the aggregate utility of their utility functions and the utility functions of their children through Bayesian optimization.
In this section it will be shown that given the global optima of the utility functions, \DBay will find the global optimum of the objective function.

\structure{Leaf agents will produce optimal assignments and utility values given a sample}
Within the sample phase of \DBay, none of the agents can optimize their variables without interaction with other agents.
The interaction involves the sending of (top-down) \sampleMessageText messages and (bottom-up) \utilityMessageText messages.
Therefore, for the leaf agents, the optimization depends on their utility functions and the \sampleMessageText message of their parent, \sampleMessageIndexedParent, as
\begin{align}
    \label{eq:optimal_leaf_agent_assignment}
    \hat{\assignment}_{\localVariablesIndexed}
    =
    \displaystyle\argmin_{\assignment \in \searchSpace_{\localVariablesIndexed}}
    \displaystyle\operator
    \Big(
        \optimumMessageLocal
    \Big)
    =
    \displaystyle\argmin_{\assignment \in \searchSpace_{\localVariablesIndexed}}
    \displaystyle\operator
    \left(
        \displaystyle\operator_{\utilityFunctionIndexed \in \functionSetLocalCombinedIndexed}
        \Big(
            \utilityFunctionIndexed(\assignmentPartial
            ~|~
            \sampleMessageIndexedParent
        \Big)
    \right)
    \qquad \forall \agentIndexed ~:~ \ChildIndexed = \emptyset.
\end{align}

When the kernel and acquisition function are selected as detailed in \cref{subsection:convergenceofBO}, the assignment $\hat{\assignment}_{\localVariablesIndexed}$ is optimal with respect to the assignments of the (pseudo-)parents of agent \agentIndexed, since $\sampleMessageIndexedParent = \hat{\assignment}_{\ParentIndexed} \cup \hat{\assignment}_{\PseudoParentIndexed}$.
Consequently, the optimal assignment results in the optimal value for the \utilityMessageText message \optimumMessageToParent given the \sampleMessageText message \sampleMessageIndexedParent.

\structure{All other agents therefore receive optimal utility values}
This optimal value is sent as a \utilityMessageText message to the parents of the leaf agents and results in the following assignment for the other agents:
\newcommand{\optimumMessageBestOptimal}{\mathcmd{
  \hat{\optimumMessage}^*
}}
\begin{align}
    \label{eq:optimal_agent_assignment}
    \hat{\assignment}_{\localVariablesIndexed}
    =
    \displaystyle\argmin_{\assignment \in \searchSpace_{\localVariablesIndexed}}
    \displaystyle\operator
    \Big(
        \optimumMessageLocal,
        \optimumMessageBestChild
    \Big)
    =
    \displaystyle\argmin_{\assignment \in \searchSpace_{\localVariablesIndexed}}
    \displaystyle\operator
    \left(
        \displaystyle\operator_{\utilityFunctionIndexed \in \functionSetLocalCombinedIndexed}
        \Big(
            \utilityFunctionIndexed(\assignmentPartial
            ~|~
            \sampleMessageIndexedParent
        \Big)
        ,
        \optimumMessageBest_\agentIndex
    \right)
    \qquad \forall \agentIndexed ~:~ \ChildIndexed \neq \emptyset.
\end{align}
The aggregated \utilityMessageText message $\optimumMessageBest_\agentIndex = \displaystyle\operator_{\agent_\agentIndexChild \in \ChildIndexed} \left( \optimumMessageFromChild \right)$ is the aggregate of the optimal \utilityMessageText messages of all children given an assignment of \agentIndexed.
Therefore, agent \agentIndexed is able to calculate the optimal assignment with respect to its parent \sampleMessageText message.

\structure{Finally, the root agent optimizes based on all aggregate \utilityMessageText messages}
This process is repeated until the root agent $\agent_\agentIndexRoot$ has received all \utilityMessageText messages from its children.
Since the root agent does not have any (pseudo-)parents, \cref{eq:optimal_agent_assignment} can be rewritten as
\begin{align}
    \hat{\assignment}_{\localVariables_\agentIndexRoot}
    =
    \displaystyle\argmin_{\assignment \in \searchSpace_{\localVariables_\agentIndexRoot}}
    \displaystyle\operator
    \Big(
        \optimumMessage_\agentIndexRoot,
        \optimumMessageBest_\agentIndexRoot
    \Big)
    =
    \displaystyle\argmin_{\assignment \in \searchSpace_{\localVariables_\agentIndexRoot}}
    \displaystyle\operator
    \left(
        \displaystyle\operator_{\utilityFunctionIndexed \in \functionSetLocalCombinedRoot}
        \Big(
            \utilityFunctionIndexed(\assignmentPartial)
        \Big)
        ,
        \optimumMessageBest_\agentIndexRoot
    \right)
    =
    \assignment^*_{\localVariables_\agentIndexRoot}.
\end{align}
Note that $\optimumMessageBest_\agentIndexRoot$ holds the aggregate utility value of all other agents based on the sample of the root agent.
For that reason, if the root agent finds the optimal assignment $\hat{\assignment}_{\localVariables_\agentIndexRoot}$ it is equal to the optimum of the objective function $\assignment^*_{\localVariables_\agentIndexRoot}$.

After the root agent has found the optimal assignment of its variables it starts the final phase of \DBay.
In this phase the root agent sends the optimal assignment in a \sampleFinalMessageText message to its children, $\sampleMessageBestRoot = \{ \assignment^*_{\localVariables_\agentIndexRoot} \}$.
Based on that optimal sample all agents are able to determine their optimal assignments, as shown in \cref{eq:optimal_agent_assignment}, and append their optimal assignment to the final message before sending the final message to their children, i.e. $\sampleMessageBestIndexed = \{ \assignment^*_{\localVariablesIndexed} \} \cup \sampleMessageBestParent$.
This process is repeated until the leaf agents are reached and all agents have assigned the global optimal values to their variables.

\subsection{Summary}
\structure{conclusion for a f with L we can select the kernel with parameter to guarantee global convergence!}
In \cref{subsection:convergenceofBO} it was shown that, based on the Lipschitz constant of a utility function (or aggregate of functions), the kernel and acquisition function (and their parameters) can be appropriately selected to guarantee convergence to the global optimum of the utility function.
Subsequently, \cref{subsection:globalconvergenceofDBAY} has shown that, if the agents are able to find the global optima of the aggregate of the utility functions in finite time, \DBay will converge to the global optimum of the objective function.
Combining these results proves the convergence of \DBay to the global optimum of the objective function for utility functions with known Lipschitz constants.

\section{Simulation Results}
\label{section:simulation_results}
\newcommand{\sensorViewAngle}{\mathcmd{\beta}}
\newcommand{\sensorRangeDistance}{\mathcmd{l}}
\newcommand{\sensorOrientation}{\mathcmd{\omega}}

In this section, the performance of \DBay is empirically evaluated using two metrics.
The first is the achieved relative utility as a function of the number of samples.
The number of samples is used as the threshold for the agents in the sample phase.
The relative utility is defined as the achieved utility divided by the optimal utility.
The relative utility allows for the comparison of the results over various randomly generated problems.
The second metric focuses on the sample efficiency of the algorithm.
This metric is important to consider if the evaluation of the utility functions is computationally expensive.
This metric compares the number of samples required by \DBay and the number of samples required by the centralized approach to achieve equal relative utility.
In the centralized approach, the continuous domains are equidistantly discretized.
The optimum of all possible combinations of these discrete domains is then taken as the achieved centralized utility.

\subsection{Sensor Coordination Problem}
The performance of \DBay is evaluated based on the sensor coordination problem.
The sensor coordination problem is an optimization problem in which every agent needs to orient its sensor in order to observe targets as accurately as possible.
This problem is modeled within the DCOP framework as a distributed problem.
A real-world analogue would be the optimization of the orientation of multiple cameras based on image recognition.
The image recognition process can require significant computational effort depending on the image quality and the type of target.
This makes it computationally intensive to check every orientation of the camera image, especially for a centralized approach.

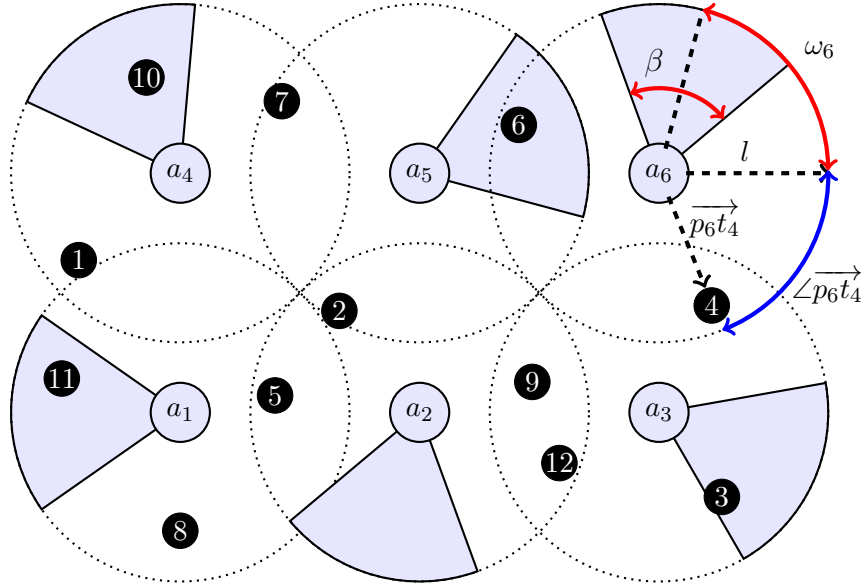
\begin{figure}[h]
	\centering
		\def\centerarc[#1](#2)(#3:#4:#5)
    { \draw[#1] ($(#2)+({#5*cos(#3)},{#5*sin(#3)})$) arc (#3:#4:#5); }

\def\centerarcslice[#1](#2)(#3:#4:#5)
    { \draw[#1] (#2) -- ($(#2)+({#5*cos(#3)},{#5*sin(#3)})$) arc (#3:#4:#5) -- (#2); }

\def\sensorRange{2.5cm}
\def\agentDistance{{1.4142 * 2.5cm}}
\def\agentDistanceTwo{{2 * 1.4142 * 2.5cm}}
\def\agentDistanceTotal{{1.4142 * 2.5cm + \sensorRange}}
\def\sensorView{35}

\def\drawAgent(#1)(#2)(#3)(#4)(#5)(#6)
    {
    \begin{scope}[xshift=#1, yshift=#2]
        \centerarcslice[thick, -, fill=blue!10](0,0)(#5-#4:#5+#4:#3);

        \centerarc[thick, dotted](0,0)(0:360:#3);

        \node[circle, minimum size=0.75cm, thick, draw, fill=blue!10] () at (0,0) {$\agent_{#6}$};

        \ifthenelse{\equal{#6}{6}}{
            \centerarc[red, ultra thick, <->](0,0)(#5+#4:#5-#4:#3/2);
            \node[] () at ({.65 * #3 * cos(#5+#4/2)},{.65 * #3 * sin(#5+#4/2)}) {\large$\sensorViewAngle$};

            \draw[ultra thick, dashed, black] ({0.15 * #3 * cos(#5)},{0.15 * #3 * sin(#5)}) -- ({#3 * cos(#5)},{#3 * sin(#5)});
            \draw[ultra thick, dashed, black, ->] (0.4,0) -- ({(0.35 + #3+1)/2},0) node[above, black] {\sensorRangeDistance} -- (#3+1,0);

            \draw[ultra thick, dashed, black, ->] ({0.15 * #3 * cos(-68)},{0.15 * #3 * sin(-68)}) -- ({0.3 * #3 * cos(-68)},{0.3 * #3 * sin(-68)}) node[right, black] {$\overrightarrow{p_{6}t_{4}}$} -- ({0.73 * #3 * cos(-68)},{0.73 * #3 * sin(-68)});
            \node[] () at ({1.2 * #3 * cos(-68 / 2)},{1.2 * #3 * sin(-68 / 2)}) {$\angle \overrightarrow{p_{6}t_{4}}$};
            \centerarc[blue, ultra thick, <->](0,0)(1:-68:#3);

            \node[] () at ({1.2 * #3 * cos(#5 / 2)},{1.2 * #3 * sin(#5 / 2)}) {$\sensorOrientation_{#6}$};
            \centerarc[red, ultra thick, <->](0,0)(1:#5:#3);
        }{}
    \end{scope}
    }

\def\drawTarget(#1)(#2)(#3)
        {
        \begin{scope}[xshift=#1, yshift=#2]
            \node[circle, scale=1.3, draw, fill=black, text=white] () at (0,0) {};
            \node[scale=1, text=white] () at (0,0) {$#3$};
        \end{scope}
        }

%
%
\begin{tikzpicture}[auto, scale=0.9]


    \drawAgent(0)(0)(\sensorRange)(\sensorView)(-180)(1)
    \drawAgent(\agentDistance)(0)(\sensorRange)(\sensorView)(-105)(2)
    \drawAgent(\agentDistanceTwo)(0)(\sensorRange)(\sensorView)(-25)(3)

    \drawAgent(0)(\agentDistance)(\sensorRange)(\sensorView)(120)(4)
    \drawAgent(\agentDistance)(\agentDistance)(\sensorRange)(\sensorView)(20)(5)
    \drawAgent(\agentDistanceTwo)(\agentDistance)(\sensorRange)(\sensorView)(75)(6)

    \drawTarget(-1.75cm)(0.5cm)(11)
    \drawTarget(-1.5cm)(2.25cm)(1)
    \drawTarget(-0.5cm)(4.95cm)(10)
    \drawTarget(0.0cm)(-1.75cm)(8)
    \drawTarget(1.4cm)(0.25cm)(5)
    \drawTarget(1.5cm)(4.6cm)(7)
    \drawTarget(2.35cm)(1.5cm)(2)
    \drawTarget(5.6cm)(-0.75cm)(12)
    \drawTarget(5cm)(4.25cm)(6)
    \drawTarget(5.2cm)(0.45cm)(9)
    \drawTarget(7.85cm)(1.58cm)(4)
    \drawTarget(8cm)(-1.25cm)(3)

\end{tikzpicture}
	\caption{
        Graphical example of a sensor coordination problem with \num{6} sensors and \num{12} targets.
        The sensors $\agentIndexed$ are arranged in a grid.
        The distance between the sensors is based on the sensor range \sensorRangeDistance.
        The sensor range is indicated by a dotted circle centered around the position of the sensor.
        The observed area of the sensors is shown as shaded areas, and is based on the angle of view \sensorViewAngle and the orientation $\sensorOrientation_\agentIndex$.
        The targets are denoted as annotated black circles.
    }
    \label{fig:sensor_coordination_problem_example}
\end{figure}

Within the sensor coordination problem, all sensors are identical in terms of their sensor range $\sensorRangeDistance$ and angle of view $\sensorViewAngle$.
These properties, combined with the position of the sensor, determine the observation domain of the sensor.
This domain defines all locations that could be observed by the sensor.
The orientation $\sensorOrientation_\variableIndex$ of the sensor \variableIndex changes the observed area within the observation domain.
A target is detected when it is located within this area.
For every detected target, positive utility is allocated to the agent responsible for the sensor.
The maximum utility is allocated when the sensor is oriented directly at the target.
The utility value decreases linearly towards the edges of the observation area.
The optimal utility is determined by the optimal solution of a centralized approach with \mbox{\num{720} samples} for every domain.
The parameters of the problem are the number of targets, the number of sensors, the sensor range \sensorRangeDistance, and the angle of view of the sensors \sensorViewAngle.
The sensors are arranged in an equally distanced rectangular grid.
The distance between the sensors is such that the combined observation domains of all sensors is maximized without allowing unobservable areas between the sensors.
While the locations of the targets $t$ are randomly distributed within the combined observation domains of the sensors.
In the simulations, all problems are generated with \num{6} sensors, \num{12} targets, and with identical sensor properties where the sensor range is set to $\sensorRangeDistance = 1$ and the angle of view is set to $\sensorViewAngle = \SI{36}{\degree}$.
A graphical example of the simulated sensor coordination problems can be seen in \cref{fig:sensor_coordination_problem_example}.

The sensor coordination problem is described within the continuous DCOP framework as follows:
\begin{itemize}
    \item \agentSetLong is the set of agents,
        where \numberOfAgents is the number of agents.
        The position of an agent is denoted as $p_\agentIndex \in \realsSet^2$.

    \item $\variableSet = \{ \sensorOrientation_1, \dots, \sensorOrientation_\numberOfVariables \}$ is the set of sensor orientations, where $\numberOfVariables = \numberOfAgents$.

    \item \domainSetLong, where $\domainSet_\agentIndex = (\SI{-180}{\degree}, \SI{180}{\degree})$ for all $\variableIndex = 1, \dots, \numberOfVariables$ indicating all possible orientations of the sensor.

    \item $\functionSet = \{f_{\functionIndex}\}_{\functionIndex=1}^{\numberOfFunctions}$ is the set of utility functions associated with the observation of the targets.
    The number of targets is denoted by $T \in \mathds{N}$.
    A target is located at position $t_\functionIndex \in \realsSet^2$.
    The utility functions of the targets are described as,
    $f_{\functionIndex} = \max_{\variableIndex = 1, \dots, \numberOfVariables}\left( f_{\functionIndex,\variableIndex} \right)$ for $\functionIndex = 1, \dots, \numberOfFunctions$, where
    \begin{align*}
        f_{\functionIndex,\agentIndex} =
        \begin{cases}
            1 - | \sensorOrientation_\variableIndex - \angle \overrightarrow{p_{\agentIndex}t_{\functionIndex}} | / \sensorViewAngle
            &\mbox{if } \| \overrightarrow{p_{\agentIndex}t_{\functionIndex}} \| \leq \sensorRangeDistance \mbox{ and } | \sensorOrientation_\variableIndex - \angle \overrightarrow{p_{\agentIndex}t_{\functionIndex}} | \leq \sensorViewAngle
            \\
            0 &\mbox{otherwise}
        \end{cases}
    \end{align*}
    and $\overrightarrow{p_{\agentIndex}t_{\functionIndex}}$ denotes the vector between the location of the target $t_{\functionIndex}$ and the position of the agent \agentIndexed.

    \item $\variableMapping(\sensorOrientation_\variableIndex) = \agentIndexed$ for $\variableIndex = 1, \dots, \numberOfVariables$ allocating a single sensor to every agent.

    \item $\operator = \sum(\cdot)$, resulting in the goal function $G(\cdot) = \sum_{f_\functionIndex \in F} f_\functionIndex(\cdot)$.

\end{itemize}

\subsection{Results}
The performance results of \DBay compared to the centralized approach are presented in \cref{fig:simulation_results}.
This figure shows the average results for 30 randomly generated problems for 6 sensors and 12 targets.
The results show an increase in the achieved relative utility of \DBay compared to the centralized approach based on the number of samples.
The difference in achieved utility can be explained by investigating the sampling strategies.
The centralized approach samples the sensor orientations equidistantly.
Therefore, as the number of samples is increased, the resolution of the samples decreases uniformly for the centralized approach.
While \DBay samples dynamically to balance exploration and exploitation based on all previously acquired observations.
Consequently, \DBay will initially focus on exploration and eventually focus on exploitation.
This behavior is clearly visible in \cref{fig:simulation_results_utility} in the range between \num{3} and \num{10} samples.
Within this range \DBay samples the sensor orientations equidistantly focussing on exploration.
This sampling behavior is identical to the centralized approach, which can be seen in the similarity in achieved utility.
At \num{11} number of samples the achieved utility of \DBay increases substantially.
This can be explained based on the angle of view of \SI{36}{\degree} of the sensors during the experiments.
At \num{10} samples the entire observation domain of a sensor was observed.
This enabled the switch to exploitation of the observations during the sample generation.
It clearly shows the advantage of the dynamic sampling of \DBay over equidistantly sampling.
The advantage is even more prominent when comparing the required number of samples by the centralized approach to achieve equal utility to \DBay, as shown in \cref{fig:simulation_results_samples}.

\newcommand{\figResultHeigth}{7cm}
\begin{figure}[h]
    \centering
    \begin{subfigure}[t]{0.5\textwidth}
		\centering
		\includegraphics[height=\figResultHeigth]{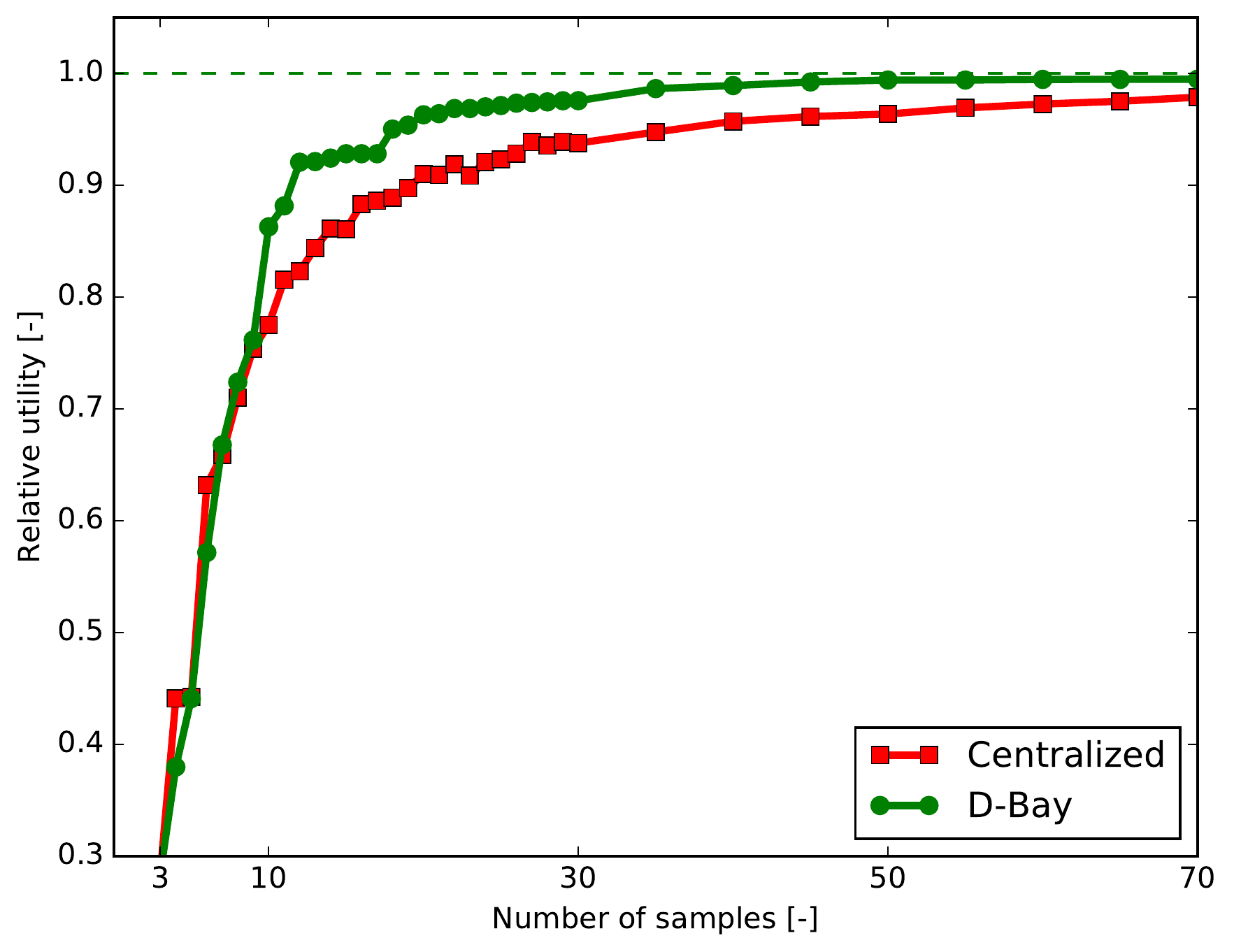}
        \caption{Achieved relative utility}
        \label{fig:simulation_results_utility}
	\end{subfigure}
	\begin{subfigure}[t]{0.45\textwidth}
		\centering
		\includegraphics[height=\figResultHeigth]{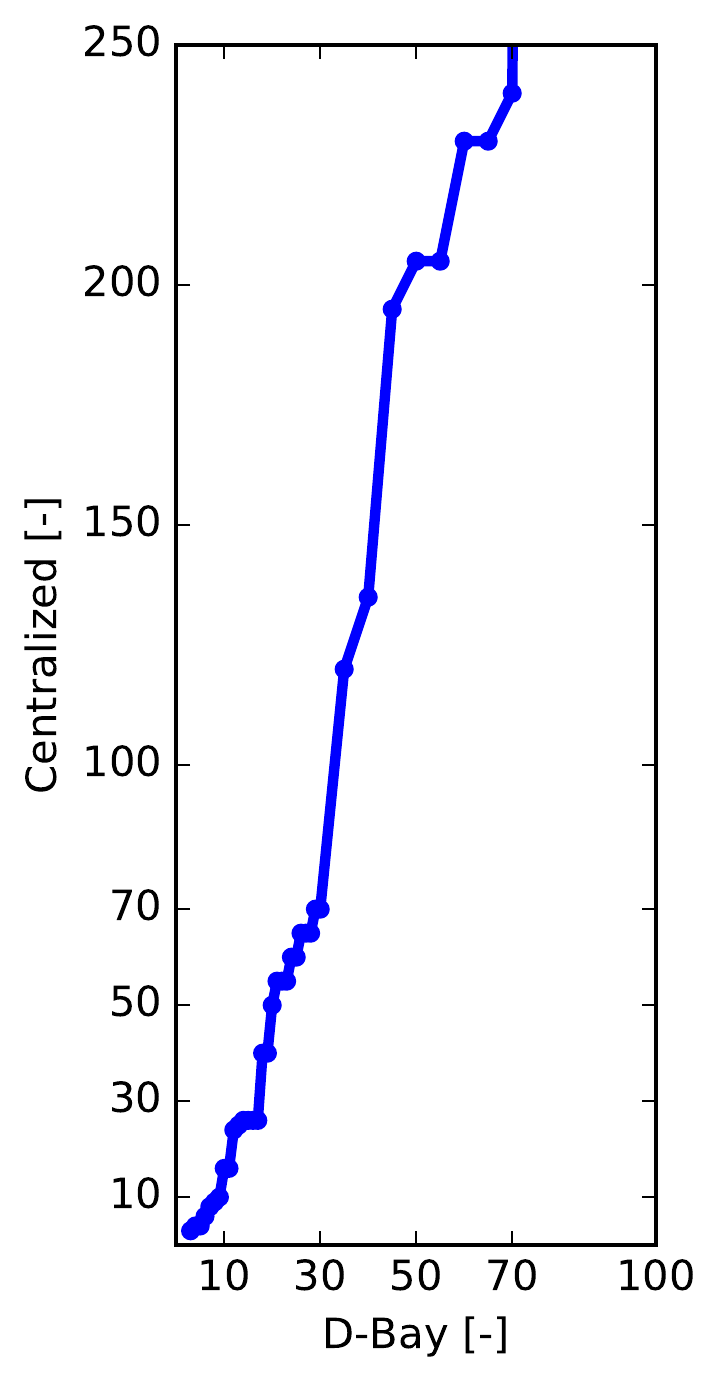}
        \caption{Required number of samples}
        \label{fig:simulation_results_samples}
	\end{subfigure}
    \caption{
        Simulation results for randomly generated sensor coordination problems with \num{6} sensors and \num{12} targets.
        The figures show the average result of \num{30} randomly generated problems.
        \cref{fig:simulation_results_utility} shows the achieved utility of \DBay and the centralized approach relative to the optimum.
        Note that the achieved utility for both algorithms is equal for \num{3} samples, since the first \num{3} samples are equidistantly spaced for \DBay.
        \cref{fig:simulation_results_samples} shows the number of samples required by the centralized approach to achieve the same utility as \DBay.
    }
    \label{fig:simulation_results}
\end{figure}

\clearpage
\section{Conclusions}
\label{section:conclusion}
In this work, the novel algorithm called \DBayLong has been introduced to solve a \DCOPLong with continuous domains.
Within \DBay, the continuous domains are sampled based on \BO.
This removes the need for discretization of the domains.
Compared to traditional DCOP solvers, which require discretization, it results in a reduction in the computational demands of the individual agents.
For utility functions with known Lipschitz constants, \DBay is proven to converge to the global optimum solution of the DCOP.

A sensor coordination problem has been used to evaluate the performance of \DBay.
The results show that \DBay outperforms a centralized approach based on the achieved utility as a function of the required number of samples.
This sample efficiency is a result of the application of \BO with \DBay.

In future work, \DBay will be extended towards dynamic DCOPs \cite{fioretto2019distributedreview} in which the agents need to optimize a dynamic problem at every time step.
An extension will increase the applicability of the developed algorithm to dynamic real-world problems in which tracking of targets is an important factor, such as multi-agent surveillance.
In a dynamic adaptation of the sensor coordination problem, the locations of the targets change over time based on the target properties, such as speed and turn radius.
Furthermore, an extension towards non-ideal communication will be investigated.
Within that extension, the (pseudo-tree) communication structure will become dynamically generated and the message passing between the agents will be asynchronous.
An additional extension would be to further decrease the computational complexity of the algorithm by removing the optimization of the acquisition function.
A possible direction would be to apply a a bound-based method such as the approach of \citeA{Kawaguchi2016}.

\section*{Acknowledgment}
The authors would like to thank dr. Peyman Mohajerin Esfahani of the Delft Center for Systems and Control at the Delft University of Technology for his insightful comments.

\clearpage  
\appendix
\section{Distributed Bayesian Algorithm}
\label{section:distributed_bayesian_algorithm}

\newcommand{\algrule}[1][.2pt]{\par\vskip.5\baselineskip\hrule height #1\par\vskip.5\baselineskip}

\begin{algorithm}

\SetKwInOut{Input}{Input}
\SetKwInOut{Output}{Output}
\Input{$\ParentIndexed, \PseudoParentIndexed, \ChildIndexed, \PseudoChildIndexed, \functionSetLocalIndexed, \functionSetSharedIndexed, \localVariablesIndexed, \kernel$}
\Output{$\hat{\assignment}_{\localVariablesIndexed}$}

\SetKwFor{ForEach}{foreach}{do}{}
\SetKwFor{When}{when received}{}{}
\SetKwBlock{Initialization}{Initialization}{}
\SetKwData{Root}{root agent}
\SetKwData{Leaf}{leaf agent}
\SetKwFunction{FnSample}{computeOptimalSample(\kernel)}
\SetKwFunction{FnStoreUtility}{storeUtility(\optimumMessageToParent, \sampleMessageIndexed)}

\SetKwFunction{FnRetrieveBestSample}{retrieveOptimalLocalSample(\sampleMessageBestIndexedParent)}
\SetKwFunction{FnRetrieveBestUtility}{retrieveOptimalUtility()}

\SetKwFunction{FnSendSample}{send(\agentIndexedChild, \sampleMessageIndexed)}
\SetKwFunction{FnSendSampleInitial}{send(\agentIndexedChild, \sampleMessageIndexed)}

\SetKwFunction{FnSendSampleBest}{send(\agentIndexedChild, \sampleMessageBestIndexed)}
\SetKwFunction{FnSendSampleBestInitial}{send(\agentIndexedChild, \sampleMessageBestIndexed)}

\SetKwFunction{FnSendUtility}{send(\ParentIndexed, \optimumMessageToParent)}

\SetKwFunction{SubOptimizeLocal}{optimizeLocalVariables(\sampleMessageIndexedParent)}
\SetKwFunction{SubOptimizeLocalRoot}{optimizeLocalVariables($\emptyset$)}

\SetKwFunction{SubSamplePhase}{createSampleMessage(\sampleMessageIndexedParent)}
\SetKwFunction{SubSamplePhaseRoot}{createSampleMessage($\emptyset$)}

\SetKwFunction{SubUtilityPhase}{calculateUtility(\optimumMessageBestChild)}
\SetKwFunction{SubUtilityPhaseLeaf}{calculateUtility()}

\SetKwFunction{SubFinalPhase}{processFinal(\sampleMessageBestIndexedParent)}
\SetKwFunction{SubFinalPhaseRoot}{processFinal($\emptyset$)}

\SetKwFunction{SubCalculateUtility}{calculateUtility(\sampleMessageIndexed)}
\SetKwFunction{SubGetChildUtility}{getChildUtility(\sampleMessageIndexed)}

\Initialization{
    \uIf{\Root}{
    \While{not threshold reached}{
        $\optimumMessageToParent \coloneqq \SubOptimizeLocalRoot$\;
    }
        \SubFinalPhaseRoot\;
    }
}
\When{\sampleMessageText \sampleMessageIndexedParent from parent \ParentIndexed}{

    \While{not threshold reached}{
        $\optimumMessageToParent \coloneqq \SubOptimizeLocal$\;
    }
    \FnSendUtility\;
}

\When{\sampleFinalMessageText \sampleMessageBestIndexedParent from parent \ParentIndexed}{
    \SubFinalPhase\;
}

\caption{\DBayLong for agent \agentIndexed}
\label{algorithm:DBAY}
\end{algorithm}

\begin{function}
    \SetKwFor{ForEach}{foreach}{do}{}
    \SetKwFor{When}{when received}{}{}
    \SetKwBlock{Initialization}{Initialization}{}
    \SetKwData{Root}{root agent}
    \SetKwData{Leaf}{leaf agent}
    \SetKwFunction{FnSample}{computeOptimalSample(\kernel)}
    \SetKwFunction{FnStoreUtility}{storeUtility(\optimumMessageToParent, \sampleMessageIndexed)}

    \SetKwFunction{FnRetrieveBestSample}{retrieveOptimalLocalSample(\sampleMessageBestIndexedParent)}
    \SetKwFunction{FnRetrieveBestUtility}{retrieveOptimalUtility()}

    \SetKwFunction{FnSendSample}{send(\agentIndexedChild, \sampleMessageIndexed)}
    \SetKwFunction{FnSendSampleInitial}{send(\agentIndexedChild, \sampleMessageIndexed)}

    \SetKwFunction{FnSendSampleBest}{send(\agentIndexedChild, \sampleMessageBestIndexed)}
    \SetKwFunction{FnSendSampleBestInitial}{send(\agentIndexedChild, \sampleMessageBestIndexed)}

    \SetKwFunction{FnSendUtility}{send(\ParentIndexed, \optimumMessageToParent)}
    \SetKwFunction{SubOptimizeLocal}{optimizeLocalVariables(\sampleMessageIndexedParent)}
    \SetKwFunction{SubOptimizeLocalRoot}{optimizeLocalVariables($\emptyset$)}

    \SetKwFunction{SubSamplePhase}{createSampleMessage(\sampleMessageIndexedParent)}
    \SetKwFunction{SubSamplePhaseRoot}{createSampleMessage($\emptyset$)}

    \SetKwFunction{SubUtilityPhase}{calculateUtility(\optimumMessageBestChild)}
    \SetKwFunction{SubUtilityPhaseLeaf}{calculateUtility()}

    \SetKwFunction{SubFinalPhase}{processFinal(\sampleMessageBestIndexedParent)}
    \SetKwFunction{SubFinalPhaseRoot}{processFinal($\emptyset$)}

    \SetKwFunction{SubCalculateUtility}{calculateUtility(\sampleMessageIndexed)}
    \SetKwFunction{SubGetChildUtility}{getChildUtility(\sampleMessageIndexed)}

    \SetKwProg{ProgSubOptimizeLocal}{Function}{}{}
    \ProgSubOptimizeLocal{\SubOptimizeLocal}{
    $\assignment_{\localVariablesIndexed} \coloneqq \FnSample$\;
    \sampleMessageIndexedOtherLongInAlgorithm\;
    $\optimumMessageToParent \coloneqq \SubCalculateUtility$\;
    \KwRet \optimumMessageToParent\;
    }
\end{function}

\clearpage  
\begin{function}
    \SetKwFor{ForEach}{foreach}{do}{}
    \SetKwFor{When}{when received}{}{}
    \SetKwBlock{Initialization}{Initialization}{}
    \SetKwData{Root}{root agent}
    \SetKwData{Leaf}{leaf agent}
    \SetKwFunction{FnSample}{computeOptimalSample(\kernel)}
    \SetKwFunction{FnStoreUtility}{storeUtility(\optimumMessageToParent, \sampleMessageIndexed)}

    \SetKwFunction{FnRetrieveBestSample}{retrieveOptimalLocalSample(\sampleMessageBestIndexedParent)}
    \SetKwFunction{FnRetrieveBestUtility}{retrieveOptimalUtility()}

    \SetKwFunction{FnSendSample}{send(\agentIndexedChild, \sampleMessageIndexed)}
    \SetKwFunction{FnSendSampleInitial}{send(\agentIndexedChild, \sampleMessageIndexed)}

    \SetKwFunction{FnSendSampleBest}{send(\agentIndexedChild, \sampleMessageBestIndexed)}
    \SetKwFunction{FnSendSampleBestInitial}{send(\agentIndexedChild, \sampleMessageBestIndexed)}

    \SetKwFunction{FnSendUtility}{send(\ParentIndexed, \optimumMessageToParent)}
    \SetKwFunction{SubOptimizeLocal}{optimizeLocalVariables(\sampleMessageIndexedParent)}
    \SetKwFunction{SubOptimizeLocalRoot}{optimizeLocalVariables($\emptyset$)}

    \SetKwFunction{SubSamplePhase}{createSampleMessage(\sampleMessageIndexedParent)}
    \SetKwFunction{SubSamplePhaseRoot}{createSampleMessage($\emptyset$)}

    \SetKwFunction{SubUtilityPhase}{calculateUtility(\optimumMessageBestChild)}
    \SetKwFunction{SubUtilityPhaseLeaf}{calculateUtility()}

    \SetKwFunction{SubFinalPhase}{processFinal(\sampleMessageBestIndexedParent)}
    \SetKwFunction{SubFinalPhaseRoot}{processFinal($\emptyset$)}

    \SetKwFunction{SubCalculateUtility}{calculateUtility(\sampleMessageIndexed)}
    \SetKwFunction{SubGetChildUtility}{getChildUtility(\sampleMessageIndexed)}

    \SetKwProg{ProgSubCalculateUtility}{Function}{}{}
    \ProgSubCalculateUtility{\SubCalculateUtility}{
    \optimumMessageLocalLongInAlgorithm\;
    \uIf{$\ChildIndexed \neq \emptyset$}{
        $\optimumMessageBestChild \coloneqq \SubGetChildUtility$\;
        \optimumMessageToParentLongInAlgorithm\;
    } \uElse{
        $\optimumMessageToParent \coloneqq \optimumMessageLocal$\;
    }
    \FnStoreUtility\;
        \KwRet \optimumMessageToParent\;
    }
\end{function}

\begin{function}
    \SetKwFor{ForEach}{foreach}{do}{}
    \SetKwFor{When}{when received}{}{}
    \SetKwBlock{Initialization}{Initialization}{}
    \SetKwData{Root}{root agent}
    \SetKwData{Leaf}{leaf agent}
    \SetKwFunction{FnSample}{computeOptimalSample(\kernel)}
    \SetKwFunction{FnStoreUtility}{storeUtility(\optimumMessageToParent, \sampleMessageIndexed)}

    \SetKwFunction{FnRetrieveBestSample}{retrieveOptimalLocalSample(\sampleMessageBestIndexedParent)}
    \SetKwFunction{FnRetrieveBestUtility}{retrieveOptimalUtility()}

    \SetKwFunction{FnSendSample}{send(\agentIndexedChild, \sampleMessageIndexed)}
    \SetKwFunction{FnSendSampleInitial}{send(\agentIndexedChild, \sampleMessageIndexed)}

    \SetKwFunction{FnSendSampleBest}{send(\agentIndexedChild, \sampleMessageBestIndexed)}
    \SetKwFunction{FnSendSampleBestInitial}{send(\agentIndexedChild, \sampleMessageBestIndexed)}

    \SetKwFunction{FnSendUtility}{send(\ParentIndexed, \optimumMessageToParent)}
    \SetKwFunction{SubOptimizeLocal}{optimizeLocalVariables(\sampleMessageIndexedParent)}
    \SetKwFunction{SubOptimizeLocalRoot}{optimizeLocalVariables($\emptyset$)}

    \SetKwFunction{SubSamplePhase}{createSampleMessage(\sampleMessageIndexedParent)}
    \SetKwFunction{SubSamplePhaseRoot}{createSampleMessage($\emptyset$)}

    \SetKwFunction{SubUtilityPhase}{calculateUtility(\optimumMessageBestChild)}
    \SetKwFunction{SubUtilityPhaseLeaf}{calculateUtility()}

    \SetKwFunction{SubFinalPhase}{processFinal(\sampleMessageBestIndexedParent)}
    \SetKwFunction{SubFinalPhaseRoot}{processFinal($\emptyset$)}

    \SetKwFunction{SubCalculateUtility}{calculateUtility(\sampleMessageIndexed)}
    \SetKwFunction{SubGetChildUtility}{getChildUtility(\sampleMessageIndexed)}

    \SetKwProg{ProgSubGetChildUtility}{Function}{}{}
    \ProgSubGetChildUtility{\SubGetChildUtility}{
    \ForEachSingleLine{$\agentIndexedChild \in \ChildIndexed$}{
        \FnSendSample\;
    }
    \When{\optimumMessageFromChild from all $\agentIndexedChild \in \ChildIndexed$}{
        \optimumMessageIndexedChildrenInAlgorithm\;
    }
        \KwRet \optimumMessageBestChild\;
    }
\end{function}

\begin{function}
    \SetKwFor{ForEach}{foreach}{do}{}
    \SetKwFor{When}{when received}{}{}
    \SetKwBlock{Initialization}{Initialization}{}
    \SetKwData{Root}{root agent}
    \SetKwData{Leaf}{leaf agent}
    \SetKwFunction{FnSample}{computeOptimalSample(\kernel)}
    \SetKwFunction{FnStoreUtility}{storeUtility(\optimumMessageToParent, \sampleMessageIndexed)}

    \SetKwFunction{FnRetrieveBestSample}{retrieveOptimalLocalSample(\sampleMessageBestIndexedParent)}
    \SetKwFunction{FnRetrieveBestUtility}{retrieveOptimalUtility()}

    \SetKwFunction{FnSendSample}{send(\agentIndexedChild, \sampleMessageIndexed)}
    \SetKwFunction{FnSendSampleInitial}{send(\agentIndexedChild, \sampleMessageIndexed)}

    \SetKwFunction{FnSendSampleBest}{send(\agentIndexedChild, \sampleMessageBestIndexed)}
    \SetKwFunction{FnSendSampleBestInitial}{send(\agentIndexedChild, \sampleMessageBestIndexed)}

    \SetKwFunction{FnSendUtility}{send(\ParentIndexed, \optimumMessageToParent)}
    \SetKwFunction{SubOptimizeLocal}{optimizeLocalVariables(\sampleMessageIndexedParent)}
    \SetKwFunction{SubOptimizeLocalRoot}{optimizeLocalVariables($\emptyset$)}

    \SetKwFunction{SubSamplePhase}{createSampleMessage(\sampleMessageIndexedParent)}
    \SetKwFunction{SubSamplePhaseRoot}{createSampleMessage($\emptyset$)}

    \SetKwFunction{SubUtilityPhase}{calculateUtility(\optimumMessageBestChild)}
    \SetKwFunction{SubUtilityPhaseLeaf}{calculateUtility()}

    \SetKwFunction{SubFinalPhase}{processFinal(\sampleMessageBestIndexedParent)}
    \SetKwFunction{SubFinalPhaseRoot}{processFinal($\emptyset$)}

    \SetKwFunction{SubCalculateUtility}{calculateUtility(\sampleMessageIndexed)}
    \SetKwFunction{SubGetChildUtility}{getChildUtility(\sampleMessageIndexed)}

    \SetKwProg{ProgFinalPhase}{Function}{}{}
    \ProgFinalPhase{\SubFinalPhase}{
    $\hat{\assignment}_{\localVariablesIndexed} \coloneqq \FnRetrieveBestSample$\;
    \sampleMessageBestIndexedLongInAlgorithm\;
    \ForEachSingleLine{$\agentIndexedChild \in \ChildIndexed$}{
        \FnSendSampleBest\;
    }
    }
\end{function}
\clearpage  
\section{Dirichlet kernel interval functions}
\label{section:dirichlet_interval_functions}
\structure{General Markovian property}
As shown in the work of \citeA[Theorem~2]{Ding2018}, a kernel \kernel  of the Markovian class reduces the mean function \meanFunctionObservationDot and the variance function \varianceFunctionObservationDot of the posterior on the interval between observations as given in \cref{eq:mean_function_interval,eq:variance_function_interval}, respectively.
For the Dirichlet kernel as defined by \cref{eq:dirichlet}, for a normalized domain $\functionInputBaseIndexed, \functionInputBaseIndexedOther \in [0, 1]$ and the kernel scale parameter \kernelScale, the non-zero elements of the $\kernelGramian^{-1}_\observationIndex(\observationSet)$ matrix are given by
\begin{align*}
    (\kernelGramian^{-1}_\observationIndex(\observationSet))_{\observationIndex,\observationIndex}
    =
    \begin{cases}
        \kernelScale^{-2}
        \frac{
            \functionInputBase_1
        }{
            \functionInputBase_1
            \bracketbig{
                \functionInputBase_2 - \functionInputBase_1
            }
        }
        , &\text{if } \observationIndex = 1,
        \\[12pt]
        \kernelScale^{-2}
        \frac{
            (\functionInputObservationRight - \functionInputObservationLeft)
        }{
            (\functionInputObservation - \functionInputObservationLeft)
            (\functionInputObservationRight - \functionInputObservation)
        }
        , &\text{if } \observationIndex \in \{ 2, \dots, \numberOfObservations - 1\},
        \\[12pt]
        \kernelScale^{-2}
        \frac{
            (1 - \functionInputBase_{\numberOfObservations-1})
        }{
            (1 - \functionInputBase_\numberOfObservations)
            (\functionInputBase_\numberOfObservations - \functionInputBase_{\numberOfObservations - 1})
        }
        , &\text{if } \observationIndex = \numberOfObservations,
    \end{cases}
\end{align*}
and
\begin{align*}
    (\kernelGramian_\observationIndex^{-1}(\observationSet))_{\observationIndex-1,\observationIndex}
    =
    (\kernelGramian_\observationIndex^{-1}(\observationSet))_{\observationIndex,\observationIndex-1}
    =
    \frac{
        -\kernelScale^{-2}
    }{
        (\functionInputObservation - \functionInputObservationLeft)
    }
    ,~\observationIndex = 2,\dots,\numberOfObservations.
\end{align*}

The mean function \meanFunctionObservationDot and the variance function \varianceFunctionObservationDot for the Dirichlet kernel can be rewritten accordingly as
\begin{align}
    \meanFunctionObservation &=
    \kernelBold^\transpose_\observationIndex (\functionInputBase, \observationSet)
    \kernelGramian^{-1}_\observationIndex (\observationSet)
    \functionOutput_\observationIndex (\observationSet)
    \nonumber \\
    &=
    \begin{bmatrix}
        0
        &
        \dots
        &
        0
        &
        \frac{\functionInputObservation - x}{\functionInputObservation - \functionInputObservationLeft}
        &
        \frac{x - \functionInputObservationLeft}{\functionInputObservation - \functionInputObservationLeft}
        &
        0
        &
        \dots
        &
        0
    \end{bmatrix}
    \begin{bmatrix}
        \functionOutputBase_1 \\
        \vdots \\
        \functionOutputObservationLeftLeft \\
        \functionOutputObservationLeft \\
        \functionOutputObservation \\
        \functionOutputObservationRight \\
        \vdots \\
        \functionOutputBase_\numberOfObservations
    \end{bmatrix}
    \nonumber\\
    &=
    \begin{bmatrix}
        \frac{\functionInputObservation - x}{\functionInputObservation - \functionInputObservationLeft}
        &
        \frac{x - \functionInputObservationLeft}{\functionInputObservation - \functionInputObservationLeft}
    \end{bmatrix}
    \begin{bmatrix}
        \functionOutputBase_{\observationIndex-1}
        \\
        \functionOutputBase_{\observationIndex}
    \end{bmatrix}
    \nonumber\\
    &=
    \frac{\functionOutputBase_{\observationIndex-1}(\functionInputObservation - x) + \functionOutputBase_{\observationIndex}(x - \functionInputObservationLeft)}{\functionInputObservation - \functionInputObservationLeft}
\end{align}

and
\begin{align}
    \varianceFunctionBase_\observationIndex(\functionInputBase | \observationSet) &=
    \kernel\left( \functionInputBase, \functionInputBase \right) - \kernelBold_\observationIndex^\transpose (\functionInputBase, \observationSet)
    \kernelGramian^{-1}_\observationIndex(\observationSet)
    \kernelBold_\observationIndex (\functionInputBase, \observationSet)
    \nonumber \\
    &=
    \kernelScale^2 \functionInputBase ( 1 - \functionInputBase )
    -
    \begin{bmatrix}
        0
        &
        \dots
        &
        0
        &
        \frac{\functionInputObservation - x}{\functionInputObservation - \functionInputObservationLeft}
        &
        \frac{x - \functionInputObservationLeft}{\functionInputObservation - \functionInputObservationLeft}
        &
        0
        &
        \dots
        &
        0
    \end{bmatrix}
    \begin{bmatrix}
        \kernelScale^2 \functionInputBase_1 (1 - \functionInputBase)
        \\
        \vdots
        \\
        \kernelScale^2 \functionInputObservationLeftLeft (1 - \functionInputBase)
        \\
        \kernelScale^2 \functionInputObservationLeft (1 - \functionInputBase)
        \\
        \kernelScale^2 \functionInputBase (1 - \functionInputObservation)
        \\
        \kernelScale^2 \functionInputBase (1 - \functionInputObservationRight)
        \\
        \vdots
        \\
        \kernelScale^2 \functionInputBase (1 - \functionInputBase_\numberOfObservations)
    \end{bmatrix}
    \nonumber \\
        &=
        \kernelScale^2
        \bracket{
            \functionInputBase ( 1 - \functionInputBase )
            -
            \bracketBig{
                \frac{\functionInputObservationLeft (1 - \functionInputBase) (\functionInputObservation - x)}{\functionInputObservation - \functionInputObservationLeft}
                +
                \frac{\functionInputBase (1 - \functionInputObservation)(x - \functionInputObservationLeft)}{\functionInputObservation - \functionInputObservationLeft}
            }
        }
    \nonumber \\
        &=
        \kernelScale^2
        \frac{
            -(\functionInputObservation - \functionInputBase)
            (\functionInputObservationLeft - \functionInputBase)
        }{
            \functionInputObservation - \functionInputObservationLeft
        }.
\end{align}

\clearpage
\vskip 0.2in
\bibliography{JAIR_Distributed_Bayesian_for_DCOPs}
\bibliographystyle{theapa}

\end{document}